\documentclass[sigconf]{acmart}
\def\BibTeX{{\rm B\kern-.05em{\sc i\kern-.025em b}\kern-.08emT\kern-.1667em\lower.7ex\hbox{E}\kern-.125emX}}
\renewcommand\footnotetextcopyrightpermission[1]{} 

\usepackage[noend]{algpseudocode}
\usepackage{array}
\usepackage{amsfonts}
\usepackage{graphicx}
\usepackage{bm,mathrsfs}
\usepackage{epstopdf}

\usepackage{tikz,pgfplots}
\usetikzlibrary{snakes,arrows,shapes,trees}
\usepackage{amsmath,amsthm}
\usepackage{amsopn}
\usepackage{listings}
\usepackage{adjustbox}
\usepackage{longtable}
\usepackage{multirow}
\usepackage{hyperref}
\usepackage{pgfplots}
\usepackage{pgfplotstable}
\usepackage{grffile}
\usetikzlibrary{arrows,shapes,plotmarks}
\usepgfplotslibrary{groupplots}
\usetikzlibrary{matrix}

\usepackage[font=small,belowskip=-8pt,aboveskip=2pt]{caption}
\usepackage{titlesec}
\newcommand{\cref}[2]{\hyperref[#2]{#1~\ref*{#2}}}
\newcommand{\figref}[1]{\hyperref[#1]{Fig.~\ref*{#1}}}
\newcommand{\secref}[1]{\hyperref[#1]{Sec.~\ref*{#1}}}
\newcommand{\tabref}[1]{\hyperref[#1]{Tab.~\ref*{#1}}}
\newcommand{\eqnref}[1]{\hyperref[#1]{Eq.~(\ref*{#1})}}
\hypersetup{
	colorlinks,
	linkcolor={blue!50!black},
	citecolor={blue!50!black},
	urlcolor={blue!80!black},
	anchorcolor = {blue!80!black},
	filecolor = {blue!80!black},
	menucolor = {blue!80!black},
	runcolor = {blue!80!black}
}

\pgfplotsset{compat=1.8}
\usepackage{soul}
\usepackage{rotating}
\usepackage{url}
\usepackage{algorithm}

\usepackage{enumitem}
\usepackage{subcaption}

\newcommand{\Vector}[1]{\mathbf{#1}}

\newcommand{\mvec}{\textsc{matvec}}

\newcommand{\tsort}{\textsc{TreeSort}}
\newcommand{\dsort}{\textsc{DistTreeSort}}

\newcommand{\Dendro}{\textsc{Dendro}}

\newcommand{\In}{\textsc{In }}
\newcommand{\Out}{\textsc{Out }}

\newcommand{\Frontera}{\href{https://frontera-portal.tacc.utexas.edu/}{Frontera}}
\newcommand{\petsc}{\href{https://www.mcs.anl.gov/petsc/}{PETSc}}
\pgfplotsset{
compat=1.8,
legend image code/.code={
\draw[mark repeat=2,mark phase=2]
plot coordinates {
(0cm,0cm)
(0.15cm,0cm)        
(0.3cm,0cm)         
};%
}
}

\newcommand{\norm}[1]{\left\lVert#1\right\rVert}
\definecolor{cpu3}{HTML}{F44336}
\definecolor{cpu4}{HTML}{2196F3}
\definecolor{cpu1}{HTML}{4CAF50}
\definecolor{cpu2}{HTML}{FFC107}
\definecolor{gpu3}{HTML}{EF9A9A}
\definecolor{gpu4}{HTML}{90CAF9}
\definecolor{gpu1}{HTML}{A5D6A7}
\definecolor{gpu2}{HTML}{FFE082}

\definecolor{cpu5}{HTML}{9932CC}

\definecolor{sq_b1}{RGB}{37,52,148}
\definecolor{sq_b2}{RGB}{44,127,184}
\definecolor{sq_b3}{RGB}{65,182,196}
\definecolor{sq_b4}{RGB}{127,205,187}
\definecolor{sq_b5}{RGB}{199,233,180}
\definecolor{sq_b6}{RGB}{255,255,204}

\definecolor{sq_r1}{RGB}{189,0,38}
\definecolor{sq_r2}{RGB}{240,59,32}
\definecolor{sq_r3}{RGB}{253,141,60}
\definecolor{sq_r4}{RGB}{254,178,76}
\definecolor{sq_r5}{RGB}{254,217,118}
\definecolor{sq_r6}{RGB}{255,255,178}

\definecolor{sq_g1}{RGB}{0,104,55}
\definecolor{sq_g2}{RGB}{49,163,84}
\definecolor{sq_g3}{RGB}{120,198,121}
\definecolor{sq_g4}{RGB}{173,221,142}
\definecolor{sq_g5}{RGB}{217,240,163}
\definecolor{sq_g6}{RGB}{255,255,204}

\definecolor{div_c1}{RGB}{230,171,2}
\definecolor{div_c2}{RGB}{102,166,30}
\definecolor{div_c3}{RGB}{231,41,138}
\definecolor{div_c4}{RGB}{117,112,179}
\definecolor{div_c5}{RGB}{217,95,2}
\definecolor{div_c6}{RGB}{27,158,119}
\definecolor{div_c7}{RGB}{215,48,39}

\definecolor{div_d1}{RGB}{215,25,28}
\definecolor{div_d2}{RGB}{253,174,97}
\definecolor{div_d3}{RGB}{255,255,191}
\definecolor{div_d4}{RGB}{171,217,233}
\definecolor{div_d5}{RGB}{44,123,182}
\definecolor{ao}{RGB}{0.0, 128, 0.0}
\newtheorem*{remark}{Remark}
\usepackage{dirtytalk}

\newcommand{\added}[1]{\textcolor{black}{#1}}
\usepackage[normalem]{ulem}

\acmConference[SC'21]{Supercomputing '21: The International Conference for High Performance Computing, Networking, Storage, and Analysis}{November 14--19, 2021}{St. Louis, MO}
\acmBooktitle{Supercomputing '21: The International Conference for High Performance Computing, Networking, Storage, and Analysis,
  November 14--19, 2021, St. Louis, MO}
\acmPrice{15.00}

\begin{document}

\title{Scalable adaptive PDE solvers in arbitrary domains}
\author{Kumar Saurabh\S}
\affiliation{%
  \institution{Iowa State University}
  \streetaddress{2529 Union Drive}
  \city{Ames}
  \state{Iowa}
    \country{U.S.A}
  \postcode{50011}
}

\author{Masado Ishii\S}
\affiliation{%
  \institution{University of Utah}
  \streetaddress{50 S Central Campus Drive}
  \city{Salt Lake City}
  \state{Utah}
    \country{U.S.A}
  \postcode{84112}
}

\author{Milinda Fernando}
\affiliation{%
  \institution{University of Utah}
  \streetaddress{50 S Central Campus Drive}
  \city{Salt Lake City}
  \state{Utah}
    \country{U.S.A}
  \postcode{84112}
}

\author{Boshun Gao}
\affiliation{%
  \institution{Iowa State University}
  \streetaddress{2529 Union Drive}
  \city{Ames}
  \state{Iowa}
    \country{U.S.A}
  \postcode{50011}
}

\author{Kendrick Tan}
\affiliation{%
  \institution{Iowa State University}
  \streetaddress{2529 Union Drive}
  \city{Ames}
  \state{Iowa}
    \country{U.S.A}
  \postcode{50011}
}

\author{Ming-Chen Hsu}
\affiliation{%
  \institution{Iowa State University}
  \streetaddress{2529 Union Drive}
  \city{Ames}
  \state{Iowa}
    \country{U.S.A}
  \postcode{50011}
}

\author{Adarsh Krishnamurthy}
\affiliation{%
  \institution{Iowa State University}
  \streetaddress{2529 Union Drive}
  \city{Ames}
  \state{Iowa}
    \country{U.S.A}
  \postcode{50011}
}

\author{Hari Sundar}
\affiliation{%
  \institution{University of Utah}
  \streetaddress{50 S Central Campus Drive}
  \city{Salt Lake City}
  \state{Utah}
    \country{U.S.A}
  \postcode{84112}
}

\author{Baskar Ganapathysubramanian}
\affiliation{%
  \institution{Iowa State University}
  \streetaddress{2529 Union Drive}
  \city{Ames}
  \state{Iowa}
  \country{U.S.A}
  \postcode{50011}
}
\renewcommand{\shortauthors}{Saurabh and Ishii, et al.}
\thanks{\S These authors contributed equally}
\begin{abstract}
 
Efficiently and accurately simulating partial differential equations (PDEs) in and around arbitrarily defined geometries, especially with high levels of adaptivity, has significant implications for different application domains. A key bottleneck in the above process is the fast construction of a `good' adaptively-refined mesh. In this work, we present an efficient novel octree-based adaptive discretization approach capable of carving out arbitrarily shaped void regions from the parent domain: an essential requirement for fluid simulations around complex objects. Carving out objects produces an \textit{incomplete octree}.  We develop efficient top-down and bottom-up traversal methods to perform finite element computations on \textit{incomplete} octrees. We validate the framework by (a) showing appropriate convergence analysis and (b) computing the drag coefficient for flow past a sphere for a wide range of Reynolds numbers ($\mathcal{O}(1-10^6)$) encompassing the \emph{drag crisis} regime. Finally, we deploy the framework on a realistic geometry on a current project to evaluate COVID-19 transmission risk in classrooms.
 
 
\end{abstract}
\maketitle
\keywords{Adaptive mesh refinement, Octrees, CFD, FEM}

\section{Introduction}
\label{sec:intro}
The discretization of the domain (i.e., mesh generation) is a critical aspect of numerically solving PDEs. The resolution and quality of the mesh are intimately related to the overall accuracy of PDE solvers. Even though mesh generation is a fundamental part of numerical approaches, creating quality meshes continues to be a significant bottleneck in the overall workflow. This bottleneck is exacerbated when considering adaptivity and parallel deployment and becomes exceptionally challenging in the presence of an arbitrarily shaped geometric object that has to be \textit{carved} out from the computational domain. Such challenges are particularly common in simulating the flow over external objects. Streamlining this workflow is one of the components of the NASA 2030 computational fluid dynamics (CFD) milestone towards the goal of conducting overnight large-eddy simulations (LES)~\citep{slotnick2014cfd}: \textit{"Mesh generation and adaptivity continue to be significant bottlenecks in the CFD workflow."} 

\begin{figure}
  \begin{subfigure}{0.35\linewidth}
    \includegraphics[width=0.9\linewidth]{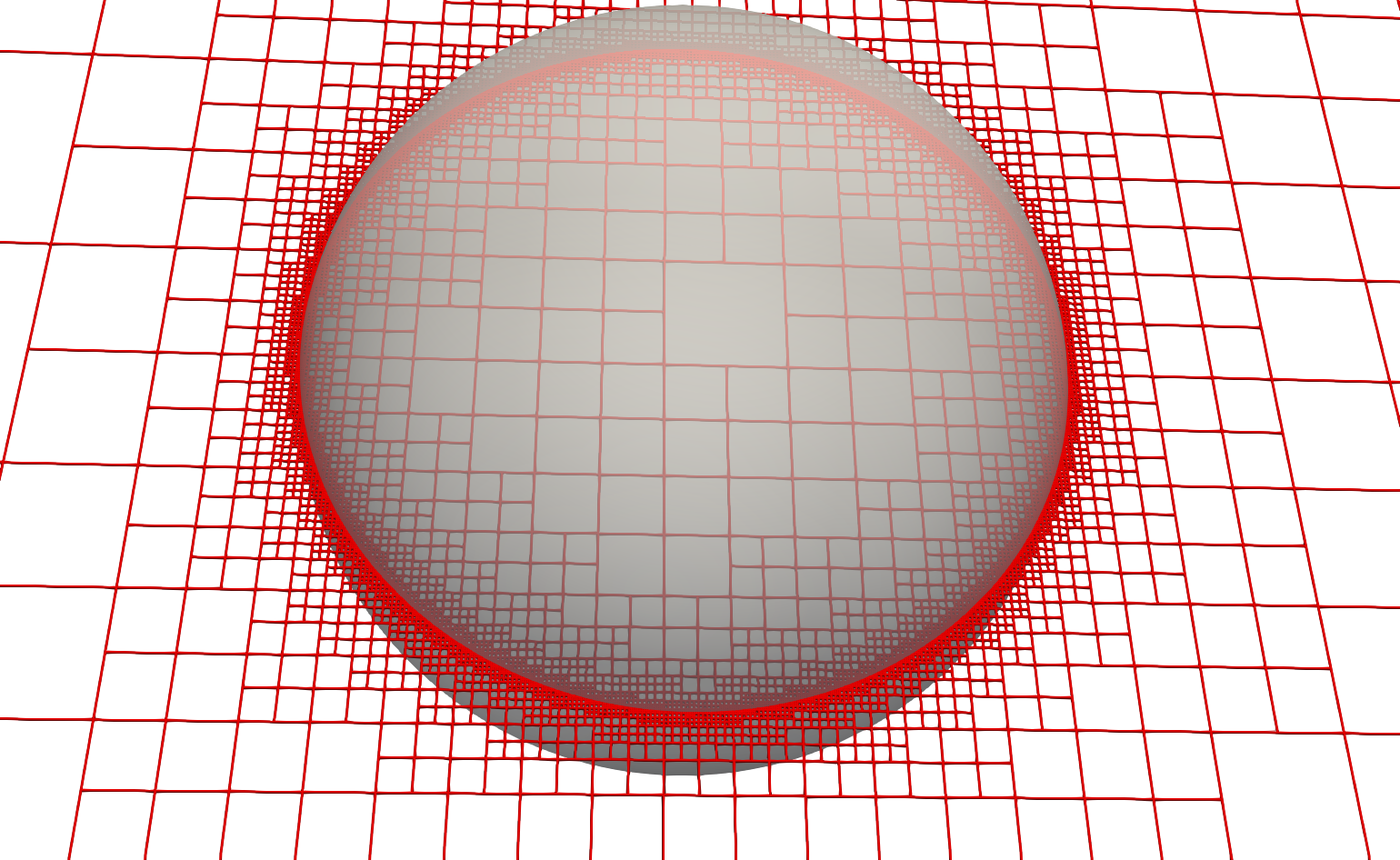}
    \caption{immersed}
    \label{fig:immersed}
\end{subfigure}
\begin{subfigure}{.35\linewidth}
  \includegraphics[width=0.9\linewidth]{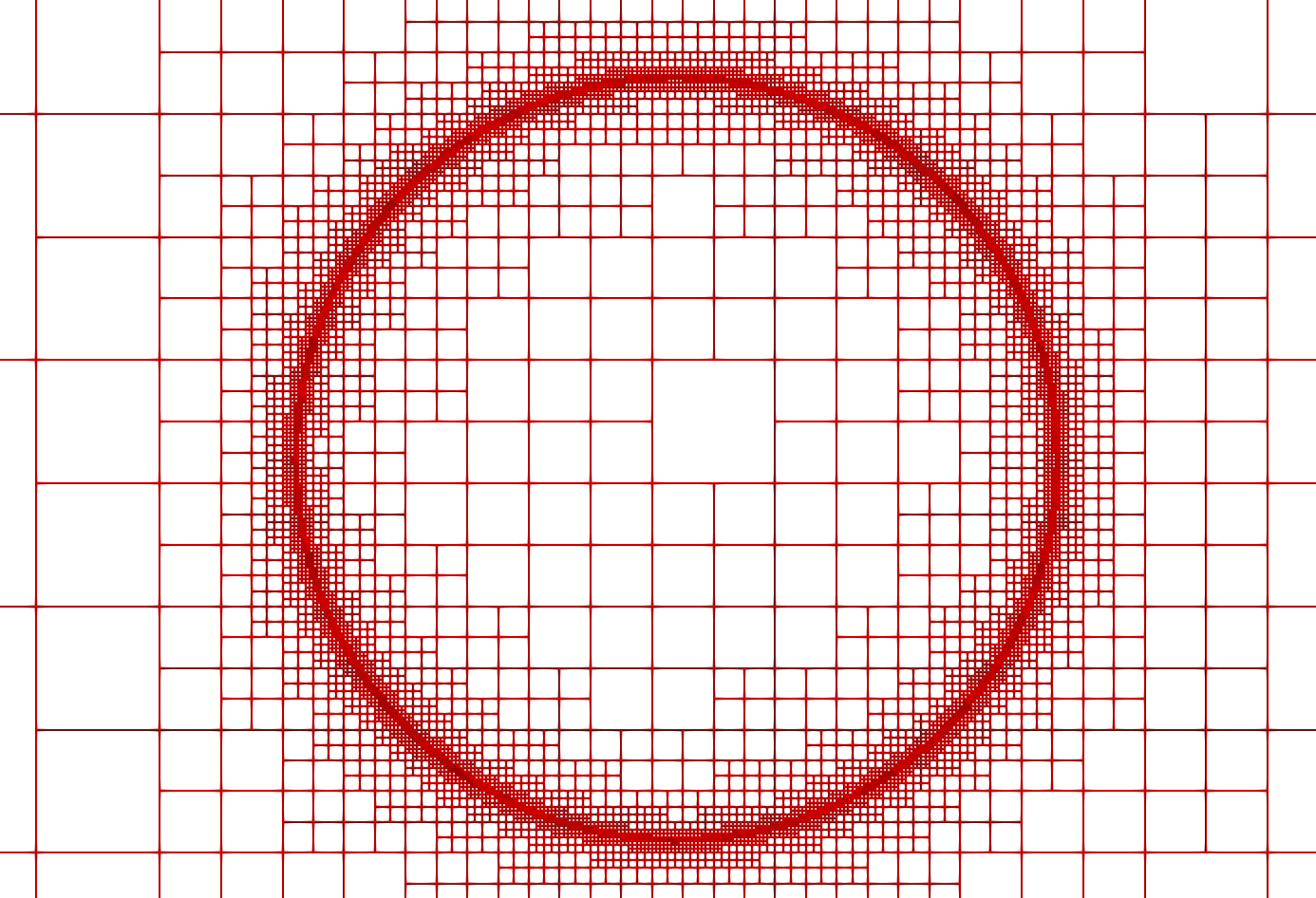}
  \caption{ complete octree}
  \label{fig:immersed-tree}
\end{subfigure}
\par\bigskip
\begin{subfigure}{0.35\linewidth}
    \includegraphics[width=0.9\linewidth]{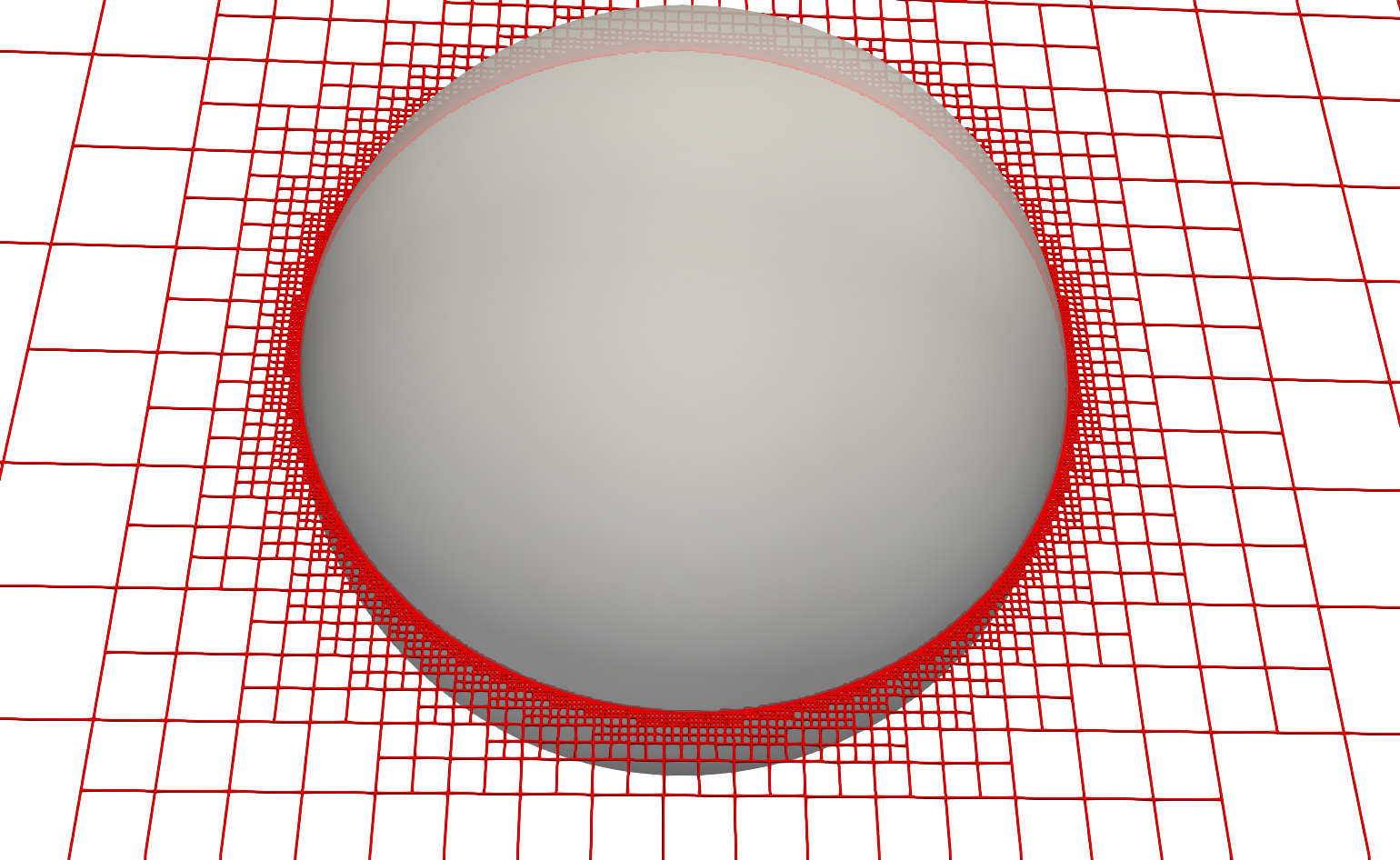}
    \caption{carved out }
    \label{fig:carved}
\end{subfigure}
\begin{subfigure}{.35\linewidth}
  \includegraphics[width=0.9\linewidth]{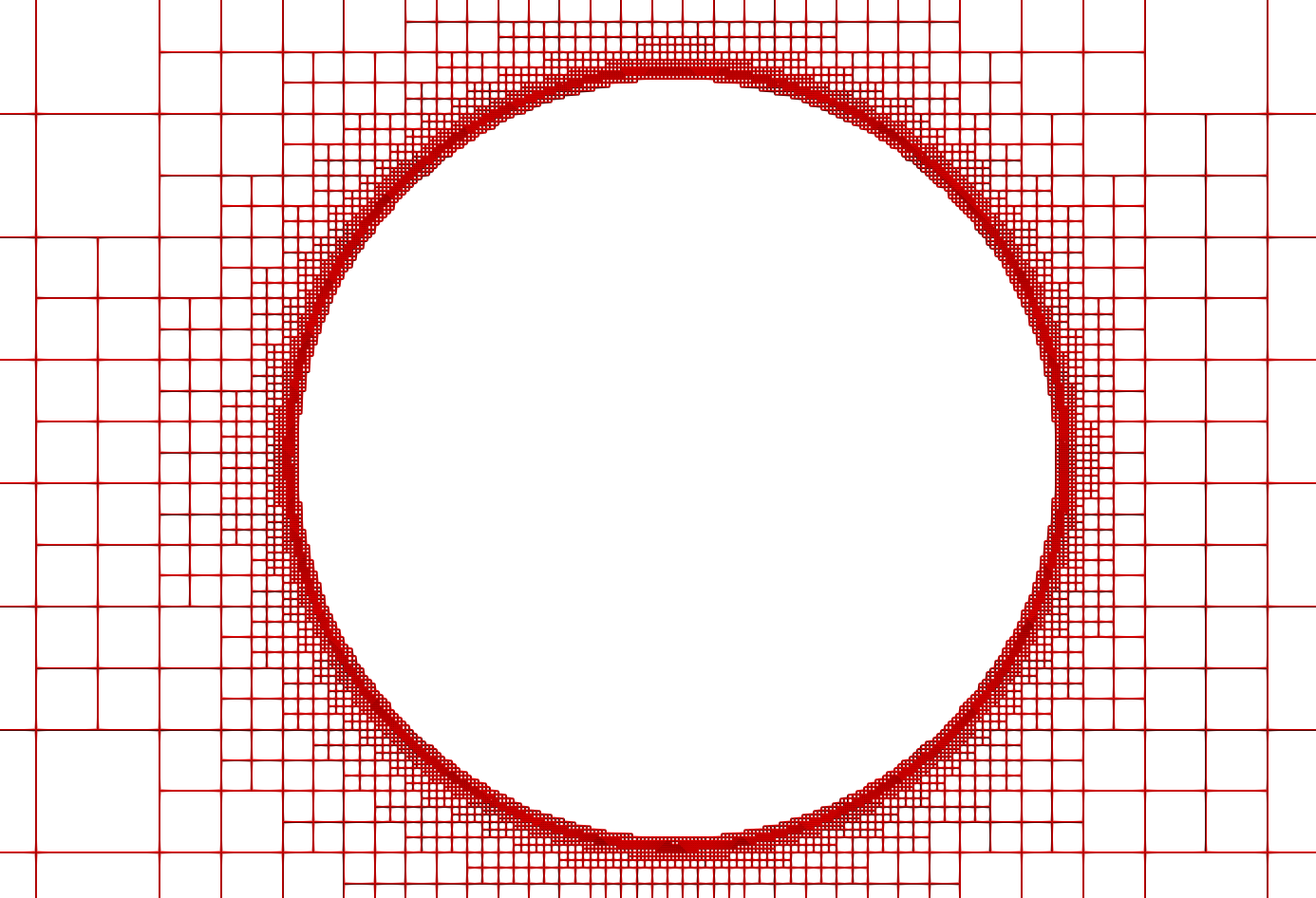}
  \caption{incomplete octree}
  \label{fig:carved-tree}
\end{subfigure}

\vspace{4 mm}
\caption{Difference between the adaptive mesh for \textit{immersed} and \textit{carved} out for the sphere case. In immersed case, we retain the full octree and this gives to a significantly large number of elements and nodes compared to the carved out case. It must be noted the elements that are inside the object do not contribute to the accuracy of the solution. Eventually Dirichlet Boundary condition are imposed on all the \In nodes. \vspace{1 mm}}
\label{fig: meshComparison}
\end{figure}

Immersed boundary methods (IBMs)~\citep{mittal2005immersed} are commonly used to simulate fluid flow around geometric objects \textit{immersed} in a computational domain. A significant advantage of IBM approaches arises from performing the complete simulation on structured grids~\citep{peskin1977numerical, mittal2005immersed}, thus avoiding any requirement of the grid conforming to the immersed geometric object (\figref{fig:immersed}). 
Naively immersing the object can lead to large void regions. These regions do not participate in the solution but require the associated matrix and vector storage as they form a part of the mesh data structure. This problem is exacerbated in the presence of multiple objects. This paper presents a strategy wherein the void regions are first \textit{carved} out from the main computational domain (\figref{fig:carved}) ---in the vein of the finite cell approach~\citep{varduhn2016tetrahedral} with the object then immersed in the carved domain. This approach reduces the number of degrees of freedom and, hence, the memory footprint associated with the void regions.


Tree-based grid generation (quadtrees in 2D and octrees in 3D) is common in computational sciences \citep{tu2005scalable,ishii2019solving,fernando2018massively,sundar2008bottom,bastian2008generic,greaves1999hierarchical,bader2012space,BursteddeWilcoxGhattas11,popinet2003gerris} largely due to its simplicity and parallel scalability. The ability to efficiently refine (and coarsen) regions of interest using tree-based data structures have made it possible to deploy them on large-scale multi-physics simulations~\citep{akhtar2013boiling,saurabh2020industrial,fernando2018massively,khanwale2020fully,rudi2015extreme,bielak2005parallel,losasso2004simulating,kim2003large,agbaglah2011parallel,losasso2004simulating}. Existing algorithms for tree-based grid generation are mainly focused on axis-aligned hierarchical splitting on isotropic domains (i.e., spheres, squares, and cubes). Such approaches cannot easily support anisotropic domains (for example, an elongated channel) or body-conforming mesh generation for a carved-out object. Standard workarounds include stretching or warping (using a coordinate transformation) of the computational domain, transforming a cubic domain into a rectangular channel. Unfortunately, these asymmetric transformations come at the cost of degradation in the overall quality of the domain discretization leading to large condition numbers in the resultant matrix (see \secref{sec:ConditionNumber}).

Our contributions in this paper are as follows: (a) we develop an efficient tree-based adaptive mesh generation framework that relaxes the requirement of the mesh to conform to isotropic domains; (b) we compare the current approach with the state-of-the-art \textit{immersed} method strategies~\citep{xu2016tetrahedral,kang2011comparative,wang2008combined,uhlmann2005immersed}; (c) we show the parallel scalability of our framework on the~\Frontera{} supercomputer up to 16K cores; (d) we deploy the framework in conjunction with a well established FEM formulation: variational multiscale (VMS) method~\citep{bazilevs2007variational} to model non-trivial applications of simulating the flow fields in classrooms to understand the risk of transmission of Coronavirus. The fast generation of quality meshes is pivotal to this application. Here, we present an octree-based mesh generation tool that provides an alternative to using two-tier meshes (HHG~\cite{hhg04}, p4est~\cite{BursteddeWilcoxGhattas11}) and is not dependent on having top-level hexahedral meshes--that can be hard to generate. In contrast, our approach works with any arbitrary user-supplied function that returns \In or \Out (of the object) for any queried point.


\section{Related Work}
\label{sec:related_work}
There have been significant algorithmic advances for the fast generation of octrees on modern supercomputers. For instance, octrees have been used for voxelization of 3D objects to accelerate ray-tracing~\citep{whang1995octree}, signed distance calculation~\citep{yu2015scalable}, compression~\citep{samet2002octree}, and fast rendering~\citep{laine2010efficient}. Building upon these successes, octrees have become one of the more common mesh generation tools for large-scale PDE simulations, with scalable and adaptive capabilities~\citep{tu2005scalable,egan428direct,rudi2015extreme,khanwale2020fully,fernando2018massively,griffith2007adaptive}.
But most work related to solving PDEs has been focused on the generation of \textit{complete} octrees in the context of PDE solvers \cite{sampath2008dendro,sundar2008bottom,fernando2017machine,fernando2018massively,ishii2019solving,weinzierl2019peano,macri2008octree,teunissen2018afivo,camata2013parallel}.
In a complete octree, every non-leaf subtree has all $2^d=8$ children,
and thus, the union of all leaf octants is a filled cube without holes. This makes simulating non-cuboid domains non-trivial, with most approaches either relying on stretching (coordinate transforms) or using a much larger bounding box.  Secondly, complex objects have to be \textit{immersed} into the octree mesh \cite{saurabh2020industrial,egan428direct,griffith2007adaptive}, rather than being \textit{carved-out}. Naively immersing the object in the octree can leave many elements that fall into the void regions. This problem is further exacerbated in the presence of multiple objects. The elements that fall into the void regions are not solved during the simulations but have an associated memory footprint. The carving of an object leads to the construction of an incomplete octree. An octree is incomplete if there are non-leaf subtrees with one or more missing children. Finally, there has been limited work in developing octree-based mesh generation to efficiently solve PDEs over complex geometries. As stated earlier, carving out affords multiple advantages (explored in this paper) at the cost of a voxelated boundary of the object. We circumvent the voxelated boundary issue via an immersed boundary (IBM) formulation on this carved-out octree. 

\added{An alternate approach has been to use two-tier meshes \cite{hhg04,BursteddeWilcoxGhattas11} that rely on having a top-level unstructured hexahedral mesh that conforms to the complex geometry and can independently refine (uniformly or adaptively) each hexahedral element of the top-level mesh. While this approach works well for simple shapes, like spheres \cite{sc08-seismic}, hex-meshing is non-trivial for more complex geometries \cite{viertel2019approach}. An alternative is then to use affine transforms within each top-level element but limits the ability to have isoparametrically refined elements. In contrast, our approach can take an arbitrary function to carve out the domain and is capable of on-the-fly refinement and coarsening that matches the arbitrary function within the refinement tolerance.  }

In this work, we address some of the key algorithmic challenges that are important to consider for carrying out efficient carving-out within the octree framework, yet missing from the existing literature:
\begin{itemize}[topsep=3pt,leftmargin=*]
    \item Careful mathematical abstraction that ensures the correctness of the generated mesh for any given arbitrary shape. Additionally, within the correctly carved out region, it is critical to correctly mark the boundary elements and nodes to solve PDEs correctly. (see \secref{sec: mathAbstraction}) 
    \item Handling of hanging nodes is critical during the carving out operations. Specifically, no hanging nodes should be present at the boundaries. If so, the parent of these nodes can lie in the inactive region, which is discarded during tree pruning. This would result in an incorrect PDE solution. (see~\secref{createNodes})
    \item It is essential that the partitioning algorithm only looks at the active region of the octree. This ensures that FEM computations are evenly distributed and hence load-balanced. The data structure used in previous literature~\cite{xu2021octree} first builds a complete octree distributed among processors before canceling out the inactive regions. This leads to the generation of complete trees with a substantial fraction of the trees in inactive regions. (see \secref{sec: octreeConstruct})
    \item Efficient pruning of trees at coarser levels is essential. Earlier approaches first generate the complete octrees before pruning. This can lead to substantial overheads for non-cube geometries, such as an elongated channel. (see \secref{sec: octreeConstruct})
\end{itemize}

\begin{figure}
\begin{subfigure}{0.4\linewidth}
\centering
    \begin{tikzpicture}[scale=0.4, transform shape]
     \fill[cpu1] (0.0,0.0) rectangle (7.5,7.5);
    \fill[cpu2] (0.75,4.5) rectangle (5.25,5.25);
    \fill[cpu2] (1.50,5.25) rectangle (4.5,6.0);
    \fill[cpu2] (0.75,0.75) rectangle (1.5,4.5);
    \fill[cpu2] (4.5,0.75) rectangle (5.25,4.5);
    \fill[cpu2] (0.0,1.5) rectangle (0.75,4.5);
    \fill[cpu2] (5.25,1.5) rectangle (6.0,4.5);
    \fill[cpu2] (1.50,0.0) rectangle (4.5,3.0);
    
    \fill[cpu4] (4.5,1.5) rectangle (1.5,4.5);
    \fill[cpu4] (0.75,2.25) rectangle (1.5,3.75);
    \fill[cpu4] (4.5,2.25) rectangle (5.25,3.75);
    \fill[cpu4] (2.25,4.5) rectangle (3.75,5.25);
    \fill[cpu4] (2.25,0.75) rectangle (3.75,1.5);
    \draw[step=1.5cm, thick] (0,0) grid (7.5,7.5);
    \draw[-,black, thick] (0,0.75) -- (6,0.75);
    \draw[-,black, thick] (0,5.25) -- (6,5.25);
    \draw[-,black, thick] (0.75,0) -- (0.75,6.0);
    \draw[-,black, thick] (5.25,0) -- (5.25,6.0);
    \draw[-,black, thick] (2.25,0) -- (2.25,1.5);
    \draw[-,black, thick] (3.75,0) -- (3.75,1.5);
    \draw[-,black, thick] (2.25,4.5) -- (2.25,6.0);
    \draw[-,black, thick] (3.75,4.5) -- (3.75,6.0);
    
    \draw[-,black, thick] (0,2.25) -- (1.5,2.25);
    \draw[-,black, thick] (0,3.75) -- (1.5,3.75);
    
    \draw[-,black, thick] (4.5,2.25) -- (6,2.25);
    \draw[-,black, thick] (4.5,3.75) -- (6,3.75);

    \draw[red,ultra thick] (3,3 ) circle (2.5cm);
   
    \foreach \x in {0,1.5,...,7.5}{
            \draw [fill=red](\x,0.0) circle (0.1cm);
            \draw [fill=red](\x,7.5) circle (0.1cm);
            \draw [fill=red](\x,6.0) circle (0.1cm);
    }
    \foreach \y in {0,1.5,...,7.5}{
            \draw [fill=red](0.0,\y) circle (0.1cm);
            \draw [fill=red](6.0,\y) circle (0.1cm);
            \draw [fill=red](7.5,\y) circle (0.1cm);
    }
    
    \foreach \x in {1.5,3.0,4.5}
    \foreach \y in {1.5,3.0,4.5}
        \draw [fill=gray](\x,\y) circle (0.1cm);
        
    \foreach \y in {0.0,0.75,...,5.9}{
        \draw [fill=red](0.0,\y) circle (0.1cm);
        \draw [fill=red](0.75,\y) circle (0.1cm);
        \draw [fill=red](5.25,\y) circle (0.1cm);
    }
    \foreach \x in {0.0,0.75,...,5.9}{
        \draw [fill=red](\x,0.0) circle (0.1cm);
        \draw [fill=red](\x,0.75) circle (0.1cm);
        \draw [fill=red](\x,5.25) circle (0.1cm);
    }
    
    \foreach \x in {2.25,3.0,...,4.2}{
        \draw [fill=gray](\x,0.75) circle (0.1cm);
        \draw [fill=gray](\x,5.25) circle (0.1cm);
    }
    \foreach \y in {2.25,3.0,...,4.2}{
        \draw [fill=gray](5.25,\y) circle (0.1cm);
        \draw [fill=gray](0.75,\y) circle (0.1cm);
    }
    
    
    
    \end{tikzpicture}
    \caption{Complete quadtree}
    \label{fig: complete}
    \end{subfigure}
    \hspace{3 mm}
    \begin{subfigure}{0.4\linewidth}
    \centering
    \begin{tikzpicture}[scale=0.4, transform shape]
     \fill[cpu1] (0.0,0.0) rectangle (7.5,7.5);
    \fill[cpu2] (0.75,4.5) rectangle (5.25,5.25);
    \fill[cpu2] (1.50,5.25) rectangle (4.5,6.0);
    \fill[cpu2] (0.75,0.75) rectangle (1.5,4.5);
    \fill[cpu2] (4.5,0.75) rectangle (5.25,4.5);
    \fill[cpu2] (0.0,1.5) rectangle (0.75,4.5);
    \fill[cpu2] (5.25,1.5) rectangle (6.0,4.5);
    \fill[cpu2] (1.50,0.0) rectangle (4.5,3.0);
    

    \draw[step=1.5cm, thick] (0,0) grid (7.5,7.5);
    \draw[-,black, thick] (0,0.75) -- (6,0.75);
    \draw[-,black, thick] (0,5.25) -- (6,5.25);
    \draw[-,black, thick] (0.75,0) -- (0.75,6.0);
    \draw[-,black, thick] (5.25,0) -- (5.25,6.0);
    \draw[-,black, thick] (2.25,0) -- (2.25,1.5);
    \draw[-,black, thick] (3.75,0) -- (3.75,1.5);
    \draw[-,black, thick] (2.25,4.5) -- (2.25,6.0);
    \draw[-,black, thick] (3.75,4.5) -- (3.75,6.0);

    \draw[-,black, thick] (0,2.25) -- (1.5,2.25);
    \draw[-,black, thick] (0,3.75) -- (1.5,3.75);
    
    \draw[-,black, thick] (4.5,2.25) -- (6,2.25);
    \draw[-,black, thick] (4.5,3.75) -- (6,3.75);

    \draw[red,ultra thick] (3,3 ) circle (2.5cm);

    \draw[-,white, ultra thick] (2.25,1.5) -- (3.75,1.5);
    \draw[-,white, ultra thick] (0.75,3.0) -- (5.25,3.0);
    \draw[-,white, ultra thick] (2.25,4.5) -- (3.75,4.5);
    \foreach \x in {0,1.5,...,7.5}{
            \draw [fill=red](\x,0.0) circle (0.1cm);
            \draw [fill=red](\x,7.5) circle (0.1cm);
            \draw [fill=red](\x,6.0) circle (0.1cm);
    }
    \foreach \y in {0,1.5,...,7.5}{
            \draw [fill=red](0.0,\y) circle (0.1cm);
            \draw [fill=red](6.0,\y) circle (0.1cm);
            \draw [fill=red](7.5,\y) circle (0.1cm);
    }

    \foreach \y in {0.0,0.75,...,5.9}{
        \draw [fill=red](0.0,\y) circle (0.1cm);
        \draw [fill=red](0.75,\y) circle (0.1cm);
        \draw [fill=red](5.25,\y) circle (0.1cm);
    }
    \foreach \x in {0.0,0.75,...,5.9}{
        \draw [fill=red](\x,0.0) circle (0.1cm);
        \draw [fill=red](\x,0.75) circle (0.1cm);
        \draw [fill=red](\x,5.25) circle (0.1cm);
    }

    
    
     \fill[white] (4.5,1.5) rectangle (1.5,4.5);
    \fill[white] (0.75,2.25) rectangle (1.5,3.75);
    \fill[white] (4.5,2.25) rectangle (5.25,3.75);
    \fill[white] (2.25,4.5) rectangle (3.75,5.25);
    \fill[white] (2.25,0.75) rectangle (3.75,1.5);
    \fill[white] (2.25,0.75) rectangle (3.75,1.5);
    \fill[white] (0.75,2.25) rectangle (5.25,3.75);

     \foreach \y in {2.25,3.0,...,4.2}{
        \draw [fill=gray](5.25,\y) circle (0.1cm);
        \draw [fill=gray](0.75,\y) circle (0.1cm);
    }
    \foreach \x in {2.25,3.0,...,4.2}{
        \draw [fill=gray](\x,0.75) circle (0.1cm);
        \draw [fill=gray](\x,5.25) circle (0.1cm);
    }
   \draw [fill=gray](1.5,1.5) circle (0.1cm);
   \draw [fill=gray](2.25,1.5) circle (0.1cm);
   \draw [fill=gray](3.75,1.5) circle (0.1cm);
   \draw [fill=gray](4.5,2.25) circle (0.1cm);
   \draw [fill=gray](5.25,2.25) circle (0.1cm);
   \draw [fill=gray](4.5,1.5) circle (0.1cm);
   \draw [fill=gray](4.5,3.75) circle (0.1cm);
   \draw [fill=gray](1.5,1.5) circle (0.1cm);
   \draw [fill=gray](1.5,3.75) circle (0.1cm);
   \draw [fill=gray](1.5,4.5) circle (0.1cm);
   \draw [fill=gray](2.25,4.5) circle (0.1cm);
   \draw [fill=gray](3.75,4.5) circle (0.1cm);
   \draw [fill=gray](4.5,4.5) circle (0.1cm);
   \draw [fill=gray](1.5,2.25) circle (0.1cm);
    \end{tikzpicture}
    \caption{Incomplete quadtree}
    \label{fig: incomplete}
    \end{subfigure}
    \vspace{4 mm}
    \caption{
    A disk (enclosed region within the red circle) immersed in a \textit{complete} (~\figref{fig: complete}) or \textit{incomplete} (~\figref{fig: incomplete}) quadtree mesh.
    Every leaf in the tree represents an element occupying a region of space,
    which is either
    completely \textit{inside} (\textcolor{cpu4}{$\blacksquare$}) the body;
    completely \textit{outside} (\textcolor{cpu1}{$\blacksquare$}) the body;
    or \textit{intercepted} (\textcolor{cpu2}{$\blacksquare$}) by the boundary (solid red circle).
    A complete tree (\figref{fig: complete}) requires all $2^d$ children
    of each non-leaf subtree to be present,
    and thus gaps are not allowed in the middle of the mesh.
    However, useful information is only contributed by elements outside (\textcolor{cpu1}{$\blacksquare$}) or intercepted (\textcolor{cpu2}{$\blacksquare$}) by the body. \vspace{1 mm}
    %
    } 
    \label{fig:carvedFig}
\end{figure}
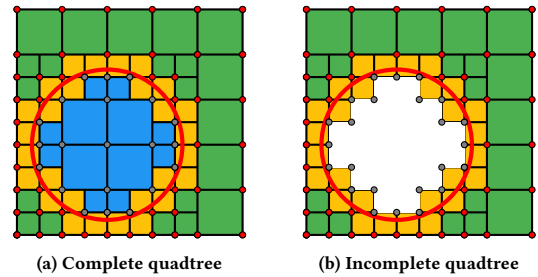{}

  

\section{Methodology}
\label{sec: Methods}
This section describes the methodology used to generate the tree-based grids (i.e., quadtrees in 2D, octrees in 3D) for arbitrary geometric domains. The key idea is to carve out regions from a $d$-dimensional cube that is \emph{inside} the immersed geometric object to generate the PDE solution domain (see \figref{fig:carvedFig}). In this paper, we refer to the aforementioned domain that is left after carving out as the \textit{subdomain}.  As noted above, the subdomain may be a sub-rectangle of a regular box, or it may have carved regions excluded from an $d$-dimensional cubic domain. The algorithms presented here are dimension agnostic, but for simplicity, we mainly focus on 3D-based grids (octree) unless specified otherwise.



Previous work~\citep{sundar2008bottom,BursteddeWilcoxGhattas11,fernando2017machine,neilsen2018massively,ishii2019solving} have demonstrated efficient methods to construct 2:1-balanced complete octrees and additional data structures to perform efficient numerical computations at a large scale. These methods order the octants of the octree using a space-filling curve (SFC) (such as the Hilbert or Morton curve) to achieve better memory accesses locality and improved distributed-memory domain decomposition. This paper presents parallel algorithms to extend numerical computations on incomplete octrees.


\subsection{Specification of the Subdomain}\label{sec: mathAbstraction}

We describe an abstraction of the application-dependent arbitrary subdomains.
The subsequent tree-based algorithms depend on a user-defined function to decide whether a given point or region in space should be retained or discarded (``carved'').
In addition, the abstraction encodes enough detail for the octree algorithms to support efficient pruning during tree traversals.

Let $\Omega = C \cup C'$ be a cube comprised of two disjoint subsets:
a closed \textit{carved} set $C \subset \Omega$, and its open complement:
the \textit{retained} set $C' \equiv \Omega \setminus C$.
Enforcing $C$ as a closed set means that it contains the boundary, $\partial C \subset C$.
(Referring back to \figref{fig:carvedFig}, in 2D,
$C$ would be the disk, including $\partial C$, the red circle.)

Suppose that $\Omega$ is hierarchically partitioned according to an octree,
$\mathcal{T}$, which captures $\partial C$ well enough under some metric.
Any octant $e$ of $\mathcal{T}$ belongs to one of the following categories,
depending on the closure of the region it bounds, $\bar{e}$:
\begin{enumerate}
    \item ``carved,'' if $\bar{e} \subset C$.
    \item ``retained,'' if $\bar{e} \not \subset C$.
    \begin{enumerate}
        \item ``intercepted,'' if ``retained'' and $\bar{e} \cap C \neq \emptyset$.
        \item ``non-intercepted,'' if $\bar{e} \subset C'$
    \end{enumerate}
\end{enumerate}

The retained octants form an incomplete octree; that is,
we define $\mathcal{T}_I \equiv \mathcal{T} \setminus \{ \text{carved leafs} \}$.
The intercepted and non-intercepted sets
specify the subdomain-boundary octants
and the subdomain-internal octants, respectively.

\subsubsection{Features of the Abstraction:\;}

Defining the subdomain abstraction in this way
ensures that the octree pruning problem is well-defined.
Notice the following:

\begin{itemize}[left=0pt,topsep=0.0in]
    \item An application can specify a subdomain through a function $F(\bar{e})$,
          where $\bar{e}$ is any filled-in cube
          of zero or positive side length, such that
    \begin{align*}
    F(\bar{e}) \coloneqq \left\{
    \begin{matrix}
      \text{label(``carved'')} && \text{if~} \bar{e} \subset $C$
      \\ 
      \text{label(``retain-internal'')} && \text{if~} \bar{e} \subset $C'$
      \\
      \text{label(``retain-boundary'')} && \text{otherwise}
    \end{matrix}
    \right.
    \end{align*}
    \begin{itemize}
        \item The function $F(\bar{e})$ applies to both octants and nodal points.
        \item Points can not be classified  ``intercepted,''
              as all points are contained in the union of $C$ and $C'$.
        \item The implementation of $F$ must take care with nontrivial intersections between $\partial C$ and an element.
        Even if all vertices lie in $C$, it is possible
        for the element to be intercepted,
        and in such a case, it should be labeled as ``retain-boundary.''
        In an application, the intersection test may be as simple  or complex
        as needed by the geometry being captured.
    \end{itemize}
        
    \item $C$ is assumed closed; hence its complement, $C^\prime$, is open.
        \begin{itemize}
            \item This convention permits the robust classification
                  of a boundary element that sits flush with $\partial C$.
                  The element is labeled ``retain-boundary,''
                  while the boundary nodes are labeled ``carved.''
            
        \end{itemize}
        \item This abstraction ensures the correctness of the generated mesh, for any arbitrary geometry, along with the correct tagging of boundary elements and nodes which is of utmost important for solving PDEs correctly.
    \item It is not necessary to generate a complete octree before filtering out the       carved octants.
          If an octant is carved, all its children are carved.
          If an octant is non-intercepted, so are all its children.
          The next section describes
          how to construct incomplete octrees
          by proactively pruning subtrees during construction.
\end{itemize}

\subsection{Octree Construction} \label{sec: octreeConstruct}

\begin{algorithm}[b!]
  \footnotesize
    \caption{ConstructUniform}
    \label{alg:construct-uniform}
    \begin{algorithmic}[1]
\Require Region $S$, SFC oracle $I$,
         final level $L$, function $F()$.
\Ensure Set $T$ of level-$L$ leafs covering subdomain, sorted by SFC.
\item[]

\If{$F(S) \neq \text{Carved}$} {\Comment Else prune}
  \If{level of $S \geq L$}
    \State $T$.push($S$)   \Comment Leaf.
    \Else
      \For{$c_\text{sfc} \leftarrow 1$ to $2^\text{dim}$}
        \Comment{Regional SFC order}
        \State $c_\text{morton} \leftarrow I$.sfc2Morton($c_\text{sfc}$)
        \State ConstructUniform($S$.child($c_\text{morton}$), $I$.child($c_\text{sfc}$))
      \EndFor
    \EndIf
\EndIf
    \end{algorithmic}
\end{algorithm}

\begin{algorithm}[t!]
 \footnotesize
    \caption{ConstructConstrained}
    \label{alg:construct-constrained}
    \begin{algorithmic}[1]
\Require Region $S$, SFC oracle $I$,
         seed octants $B$, function $F()$.
\Ensure Set $T$ of leafs, no coarser than $B$, covering the subdomain, sorted by SFC.
        
\If{$F(S) \neq \text{Carved}$} {\Comment Else prune}
  \State $L \leftarrow$ finest level in $B$
  \If{$|B| = 0$ {\bf or} level of $S \geq L$}
    \State $T$.push($S$)   \Comment Leaf.
    \Else
      \State \Comment Bucket seeds to SFC-sorted children of $S$.
      \State $l \leftarrow$ level($S$) $+1$
      \State counts[$2^\text{dim}$] $\leftarrow 0$
      \For{$b \in B$}
        \State counts[child\_num($b, l$)]$++$
      \EndFor
      \State counts[] $\leftarrow$ permute(counts, $I$)
      \State offsets[] $\leftarrow$ scan(counts)
      \State \Comment Construct child subtrees in SFC order.
      \For{$c_\text{sfc} \leftarrow 1$ to $2^\text{dim}$}
        \State $c_\text{morton} \leftarrow I$.sfc2Morton($c_\text{sfc}$)
        \State $S_c \leftarrow$ $S$.child($c_\text{morton})$
        \State $I_c \leftarrow$ $I$.child($c_\text{sfc})$
        \State $B_c \leftarrow$ $B$.slice(offsets[$c_\text{sfc}$], offsets[$c_\text{sfc}+1$])
        \State ConstructConstrained($S_c$ $I_c$, $B_c$)
      \EndFor
    \EndIf
\EndIf
    \end{algorithmic}
\end{algorithm}

\begin{algorithm}[t!]
 \footnotesize
    \caption{DistributedConstructConstrained}
    \label{alg:dist-construct-constrained}
    \begin{algorithmic}[1]
\Require Distributed set of seed octants $B$, function $F()$.
\Ensure Distributed set $T$ of leafs, no coarser than $B$, covering the subdomain, sorted by SFC.

\State \dsort($B$, load\_tol)
\Comment{\citep{ishii2019solving}}
\State $T_\text{tmp} \leftarrow$ ConstructConstrained(TreeRoot, SFC\_Root, $B$)
\State \dsort($T_\text{tmp}$, load\_tol)
\State $T_\text{local} \leftarrow$ DistributedUniqueLeafs($T_\text{tmp}$)
\State \Return $T_\text{local}$
    \end{algorithmic}
\end{algorithm}

The previous literatures have shown efficient octree algorithms to sort, construct, and traverse octrees in SFC order~\citep{sundar2008bottom,fernando2017machine,ishii2019solving}. Most of the past literature has been limited to the construction of a complete octree in an isotropic domain~\cite{sampath2008dendro,sundar2008bottom,fernando2017machine,fernando2018massively,ishii2019solving,weinzierl2019peano,macri2008octree,teunissen2018afivo,camata2013parallel}. This work builds upon them to efficiently construct the tree in presence of void regions. Algorithms~\ref{alg:construct-uniform} and~\ref{alg:construct-constrained} forms the central crux of the work, where an efficient approach for octree construction in the presence of void regions is introduced. Algorithm~\ref{alg:dist-construct-constrained} describes the partitioning algorithm based on \dsort~ to partition the trees.

The octree is constructed recursively in a top-down fashion
(see Algorithms~\ref{alg:construct-uniform} and~\ref{alg:construct-constrained}),
with child subtrees being traversed in an order determined by a regional segment of the SFC.
The top of the tree represents the entire isotropic domain, while subtrees represent cubical subregions.
A given subtree is immediately pruned if it is classified as a carved region; otherwise, it is constructed.
Constructing a subtree entails either appending the subtree as a leaf
or refining it based on a refinement criterion.
In Algorithm~\ref{alg:construct-uniform}, the refinement criterion
has a depth in the tree coarser than a target depth,
whereas in Algorithm~\ref{alg:construct-constrained}, the criterion
has a depth that is coarser than a subset of seed octants.
Other criteria are possible,
e.g., intercepting the subdomain boundary
or containing more than a maximal number of points
from an initial point cloud distribution.
If a subtree is to be refined, it is split into its eight child subtrees (in 3D).
The children are permuted into the regional SFC ordering and constructed recursively.

Algorithm~\ref{alg:dist-construct-constrained} uses the octree partitioning method based on 
\dsort~\citep{fernando2017machine,ishii2019solving}
to distribute octrees in parallel. \added{\dsort{} uses \tsort~ based comparison-free search algorithm for octree construction. Instead of performing comparison-based binary searches, \tsort~ performs MSD radix sort, except that the ordering of buckets are permuted at each level according to the specified SFC. By performing a fixed number of passes over the input data in a highly localized manner, \tsort~ avoids cache misses and random memory access  leading to  better memory performance. \cite{fernando2017machine,fernando2018massively} } Since the constructed octrees from previous algorithms entail only the active regions of the isotropic domain, \dsort~ distributes only these aforementioned active portion. This is the main difference from the past approaches, where the sorting algorithm looks at the complete tree. This step is pivotal to ensure load-balanced computation. Similar to Algorithm~\ref{alg:construct-constrained},
a set of seed octants is used to control the output tree depth.
In the distributed setting, the seed octants also inform the domain decomposition so that each rank will own approximately the same number of elements.
Note that \dsort~ accepts a tunable load-balance tolerance.
A large tolerance will partition the tree at coarse levels.
A small tolerance will balance the load more evenly
at the expense of splitting coarse subtrees over multiple processes.
Once the depth-constraining seed octants have been partitioned,
each rank constructs a tree satisfying the local constraints.
Then, overlaps between trees must be resolved.
Duplicate octants are deleted.
Finer octants are preferred to coarser overlapping octants
in order to satisfy the depth constraints globally. 

\subsection{2:1 Balancing}

\begin{algorithm}[t!]
 \footnotesize
    \caption{DistributedConstruct2to1Balanced}
    \label{alg:dist-construct-balanced}
    \begin{algorithmic}[1]
\Require Distributed set of seed octants $B$, function $F()$.
\Ensure Distributed set $T$ of leafs, no coarser than $B$, covering the subdomain, obeying 2:1-balance constraint, sorted by SFC.

\State $T_1 \leftarrow$ DistributedConstructConstrained($B$, $F$)
\State $T_2 \leftarrow$ BottomUpConstrainNeighbors($T_1$)  \Comment{$F$ not applied}
\State $T_3 \leftarrow$ DistributedConstructConstrained($T_2$, $F$)
\State \Return $T_3$
    \end{algorithmic}
\end{algorithm}

\begin{algorithm}[t!]
 \footnotesize
    \caption{BottomUpConstrainNeighbors}
    \label{alg:constrain-neighbors}
    \begin{algorithmic}[1]
\Require Unbalanced leafs $T_1$.
\Ensure Balanced seeds $T_2$.

\State $T_\text{aux}$[] $\leftarrow$ stratify $T_1$ by levels, from finest to coarsest
\For {level $l$ from finest to coarsest}
    \For {$t \in T_\text{aux}$[$l$]}
      \For {$n \in$ MakeNeighbors(MakeParent($t$))}
        \State $T_\text{aux}$[$l-1$].add\_unique($n$)
      \EndFor
    \EndFor
\EndFor
\State $T_2 \leftarrow$ concatenate($T_\text{aux}$)
\State \Return $T_2$
    \end{algorithmic}
\end{algorithm}

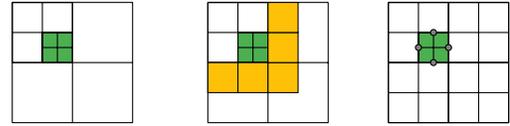
\begin{figure}
\centering
\begin{tikzpicture}[scale=0.2,every node/.style={scale=0.6} ]
\begin{scope}[shift={(25,0)}]
	\draw[step=4] (0,0) grid +(8,8);
	\draw[step=2] (0,4) grid +(4,4);
	\draw[fill=cpu1](2,4) grid +(2,2) rectangle (2,4);
\end{scope}

\begin{scope}[shift={(38,0)}]
	\draw[step=4] (0,0) grid +(8,8);
	\draw[step=2] (0,4) grid +(4,4);
	\draw[fill=cpu1] (2,4) grid +(2,2) rectangle (2,4);
	\draw[fill=cpu2] (0,2) rectangle (2,4); 
	\draw[fill=cpu2] (2,2) rectangle (4,4); 
	\draw[fill=cpu2] (4,2) rectangle (6,4); 	
	\draw[fill=cpu2] (4,4) rectangle (6,6);
	\draw[fill=cpu2] (4,6) rectangle (6,8);
\end{scope}

\begin{scope}[shift={(50,0)}]
	\draw[step=4] (0,0) grid +(8,8);
	\draw[step=2] (4,0) grid +(4,4);
	\draw[step=2] (0,4) grid +(4,4);
	\draw[step=2] (0,0) grid +(4,4);
	\draw[step=2] (4,4) grid +(4,4);
	\draw[fill=cpu1](2,4) grid +(2,2) rectangle (2,4);
	\draw [fill=gray](2,5) circle (0.2cm);
	\draw [fill=gray](4,5) circle (0.2cm);
	\draw [fill=gray](3,6) circle (0.2cm);
	\draw [fill=gray](3,4) circle (0.2cm);
\end{scope}

\end{tikzpicture}
\caption{Left most figure shows an octree which violates the 2:1 balanced constraint, where the octants that cause the violation is showed in (\textcolor{cpu1}{$\blacksquare$}). In the middle figure auxiliary balanced octants are showed in (\textcolor{cpu2}{$\blacksquare$}), in other words these are the octants needed to remove the balance constraint violation in (\textcolor{cpu1}{$\blacksquare$}). Right most figure shows the constructed octree with auxiliary balanced octants which satisfies the 2:1 balance constraint. The nodes marked by gray circles in the final 2:1 balanced mesh are hanging nodes.  \vspace{1 mm} \label{fig:aux_bal}}
\end{figure}
In a \textit{2:1-balanced} octree (\figref{fig:aux_bal}), a pair of octants sharing any parts of their boundaries
may differ in scale by at most a factor of 2:1. In other words, they may differ by at most one level in the tree.
Numerical computations on the octree grid are simpler in terms of the neighborhood data structures if the octree obeys the 2:1-balancing constraint .

We take a bottom-up approach to transform a given linear octree into a 2:1-balanced octree, based on the local block balancing method similar to the one by~\citet{sundar2008bottom}.
In our method (Algorithms~\ref{alg:dist-construct-balanced} and~\ref{alg:constrain-neighbors}),
the input octree comprises an initial set of seed octants.
The seed set is iteratively updated from the finest to the coarsest level.
For each seed octant, the neighbors of its parent octant are added to the next-coarser level of seeds.
Duplicate octants are removed from the next level before proceeding.
Finally, after all the levels have been processed, a new linear octree is constructed such that each seed octant
becomes either a leaf or an ancestor subtree in the output octree.
Thus the final seed set controls the resolution of the new octree,
ensuring the result is 2:1-balanced. It is important not to preemptively discard the carved octants,
which are generated as neighbors of parents of seed octants.
Otherwise, two leaf octants of 4:1 or greater ratio could meet in a carved region.

\subsection{Embedding Nodal Information}
\label{createNodes}
Each leaf octant in a linear 2:1-balanced octree represents an element in the FEM adaptive grid.
For a given $p$-refinement, there are $(p+1)^3$ nodes per element (in 3D).
Nodal points on the boundary of an element will be shared with same-level neighboring elements.
Nodes incident on a coarser-level neighbor is considered \textit{hanging nodes} (\figref{fig:aux_bal}).
The value of a hanging node is dependent on the values of the nodes on the coarser face or edge.
Therefore, the set of independent degrees of freedom (DOFs) on the FEM grid (underpinning a grid vector)
is defined by enumerating the unique, non-hanging nodes.

First, we loop over all elements and generate the node coordinates
with a spacing of $\left( \text{length}_\text{element} \right) / p$ in each axis.
(Nodes labeled as ``carved'' are marked as subdomain boundary nodes.)
The set of unique nodes is found by executing \tsort{~} on the nodal coordinates and removing duplicates.

An extra step is required to detect and discard hanging nodes.
In an isotropic domain, a hanging node
has fewer instances than the number expected
for an ordinary node as a function of the coordinate and grid level.
With user-specified geometry, however, the expected number of instances is nontrivial to compute.
Our solution is to explicitly ``cancel'' possible hanging nodes using temporary \textit{cancellation nodes}.
The cancellation nodes are generated on the edges and faces of elements in between the ordinary nodes,
anticipating the coordinates of hanging nodes from hypothetical finer neighbors.
After sorting, every coordinate is occupied by a mix of ordinary and cancellation nodes.
If a cancellation node is present, then the coordinate is incident on a coarser edge or face, and thus the node is hanging; the node is discarded.
Otherwise, no cancellation node is present, and the coordinate is enumerated as an ordinary node.
Thus we enumerate exactly the nodes which define a grid vector. Note that ensuring the absence of hanging nodes at the carved boundary is essential for accurate PDE solutions.

\subsection{Matrix-free, Traversal-based~\mvec}
We implement a traversal-based matrix-vector multiplication to perform matrix-free computations, which extends the methods by~\citet{ishii2019solving} to incomplete trees.

\paragraph{Matrix free:\;}
The global matrix is defined as a summation of local elemental matrices,
where the summation is due to common nodal points being shared by neighboring elements.
We are able to apply the global operator to a grid vector without explicitly assembling the global matrix.
Instead, we perform a series of elemental matrix-vector multiplications,
and use the octree structure to compose the results.

\begin{figure}[t!]
  \begin{tikzpicture}[scale=0.18,every node/.style={scale=0.6} ]
	\tikzstyle{edge from parent}=[black,->,shorten <=1pt,>=stealth',semithick,draw]
	\tikzstyle{level 1}=[sibling distance=8mm]
	\tikzstyle{level 2}=[sibling distance=6mm]
	\tikzstyle{level 3}=[sibling distance=4mm]
		
		\begin{scope}[shift={(0,0)}]
		\draw[step=10] (0,0) grid +(10,10);
		\def \r{0.12}
		\foreach \x in {0,2.5,5}{
			\foreach \y in {0,2.5,5}{
				\draw[red,fill=red] (\x,\y) circle (\r);
			}
		}
		\foreach \x in {5,7.5,10}{
			\foreach \y in {5,7.5,10}{
				\draw[red,fill=red] (\x,\y) circle (\r);
			}
		}
		\foreach \x in {0,1.25,2.5,3.75}{
			\foreach \y in {6.25,7.5,8.75,10}{
				\draw[blue,fill=blue] (\x,\y) circle (\r);
			}
		}	
		\foreach \x in {6.25,7.5,8.75,10}{
			\foreach \y in {0,1.25,2.5,3.75}{
				\draw[blue,fill=blue] (\x,\y) circle (\r);
			}
		}
		\end{scope}
		
		\begin{scope}[shift={(15,0)}]
		\draw[step=5] (0,0) grid +(10,10);
		\def \r{0.1}
		\foreach \x in {0.4,2.5,4.6,5.4}{
			\foreach \y in {0.4,2.5,4.6,5.4}{
				\draw[red,fill=red] (\x,\y) circle (\r);
			}
		}
		\foreach \x in {4.6,5.4,7.5,9.6}{
			\foreach \y in {4.6,5.4,7.5,9.6}{
				\draw[red,fill=red] (\x,\y) circle (\r);
			}
		}
		\foreach \x in {0,1.25,2.5,3.75}{
			\foreach \y in {6.25,7.5,8.75,10}{
				\draw[blue,fill=blue] (\x,\y) circle (\r);
			}
		}	
		\foreach \x in {6.25,7.5,8.75,10}{
			\foreach \y in {0,1.25,2.5,3.75}{
				\draw[blue,fill=blue] (\x,\y) circle (\r);
			}
		}
		\end{scope}
		
		\begin{scope}[shift={(30,0)}]
		\draw[step=5] (0,0) grid +(10,10);
	    \draw[step=2.5] (5,0) grid +(5,5);
		\draw[step=2.5] (0,5) grid +(5,5);
		\def \r{0.1}
		
		\foreach \x in {0.4,2.5,4.6,5.4}{
			\foreach \y in {0.4,2.5,4.6,5.4}{
				\draw[red,fill=red] (\x,\y) circle (\r);
			}
		}
		\foreach \x in {4.6,5.4,7.5,9.6}{
			\foreach \y in {4.6,5.4,7.5,9.6}{
				\draw[red,fill=red] (\x,\y) circle (\r);
			}
		}
		
		\foreach \x in {0.4,1.25,2.1,2.9,3.35,4.15}{
			\foreach \y in {6.25,7.1,7.9,8.75,9.6}{
				\draw[blue,fill=blue] (\x,\y) circle (\r);
			}
		}	
		\foreach \x in {6.25,7.1,7.9,8.75,9.6}{
			\foreach \y in {0.4,1.25,2.1,2.9,3.75}{
				\draw[blue,fill=blue] (\x,\y) circle (\r);
			}
		}
		
		\end{scope}
		\draw[->] (11,7.5) -- node[above] {top} node[below] {down} (14,7.5);
		\draw[->] (14,2.5) -- node[above] {bottom} node[below] {up} (11,2.5);	
		
		\draw[->] (26,7.5) -- node[above] {top} node[below] {down} (29,7.5);
		\draw[->] (29,2.5) -- node[above] {bottom} node[below] {up} (26,2.5);	
		
		\end{tikzpicture} 

	\caption{Illustration of top-down \& bottom-up tree traversals for a $2D$ tree with quadratic element order. The leftmost figure depicts the unique shared nodes (nodes are color-coded based on level), as we perform top-down traversal nodes shared across children of the parent get duplicated for each bucket recursively, once leaf node is reached it might be missing elemental local nodes, which can be interpolated from immediate parent (see the rightmost figure). After elemental local node computations, bottom-up traversal performed while merging the nodes duplicated in the top-down traversal.\vspace{2 mm}
	    \label{fig:traversal}
	}
\end{figure}
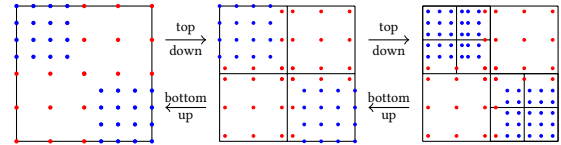

\paragraph{Traversal-based:\;}(\figref{fig:traversal})
The elemental matrix couples elemental nodes in a global input grid vector with equivalent elemental nodes in a global output grid vector.
Within a grid vector, the nodes pertaining to a particular element are generally not stored contiguously.
If one were to read and write to the elemental nodes using an element-to-node map, the memory accesses would require indirection:
$v_{glob}[map[e*npe + i]] += v_{loc}$.
Not only do element-to-node maps cause indirect memory accesses;
the maps become complicated to build if the octree is incomplete due to complex geometry.
We take an alternative approach that obviates the need for element-to-node maps.
Instead, through top-down and bottom-up traversals of the octree,
we ensure that elemental nodes are stored contiguously in a leaf,
and there apply the elemental matrix. 

The idea of the top-down phase is to selectively copy nodes
from coarser to finer levels until the leaf level,
wherein the selected nodes are exactly the elemental nodes.
Starting at the root of the tree, we have all the nodes in the grid vector.
We create buckets for all child subtrees.
Looping through the nodes, a node is copied into a bucket if the node is incident on the child subtree corresponding to that bucket.
A node that is incident on multiple child subtrees will be duplicated.
By recursing on each child subtree and its corresponding bucket of incident nodes, we eventually reach the leaf level.

Once the traversal reaches a leaf octant, the elemental nodes have been copied into a contiguous array.
The elemental matrix-vector product is computed directly, without the use of an element-to-node map.
The result is stored in a contiguous output buffer the same size as the local elemental input vector.

After all child subtrees have been traversed, the bottom-up phase returns results from a finer to a coarser level.
The parent subtree nodes are once again bucketed to child subtrees,
but instead of the parent values being copied, the values of nodes from each child are accumulated into a parent output array.
That is, for any node that is incident on multiple child subtrees, the values from all node instances are summed to a single value.
The global matrix-vector product is completed after the bottom-up phase executes at the root of the octree.

Distributed memory is supported by two slight augmentations.
Firstly, the top-down and bottom-up traversals operate on ghosted vectors.
Therefore ghost exchanges are required before and after each local traversal.
Secondly, the traversals are restricted to subtrees containing the owned octants.
The list of owned octants is bucketed top-down, in conjunction with the bucketing of nodal points.
A child subtree is traversed recursively only if one or more owned octants are bucketed to it. Note that because the traversal path is restricted by a list of existing octants,
the traversal-based~\mvec~ gracefully handles incomplete octrees without special treatment.

\begin{remark}
The traversal based~\mvec~ is designed to expose memory locality suited for deep memory hierarchies inherent in modern day clusters and accelerators like GPUs. In this work, we focus on distributed memory parallelism; the implementation on accelerators is deferred to future work.
\end{remark}

\subsection{Traversal-based Matrix Assembly}

In the previous section, we described~\mvec~ procedure that employs a tree traversal, requiring neither element-to-node maps nor global matrix assembly. In this section, we describe the matrix assembly procedure for computing the global sparse matrix. The efficient computation of matrix assembly becomes particularly important for the problems whose convergence heavily depends on the preconditioners.


To implement assembly, we have leveraged PETSc interface~\citep{petsc-web-page,petsc-efficient},
which only requires a sequence of entries
$(\text{id}_\text{row}, \text{id}_\text{col}, \text{val})$,
and can be configured to add entries with duplicate indices~\citep{petsc-user-ref}. Note that any other distributed sparse-matrix library can be supported in a similar fashion.

The remaining task is to associate the correct global node indices
with the rows and columns of every elemental matrix.
We use an octree traversal to accomplish this task.
Similar to the traversal-based~\mvec,
nodes are selectively copied from coarser to finer levels, recursively,
until reaching the leaf, wherein the elemental nodes are contiguous.
Note that integer node ids are copied instead of floating-point values from a grid vector.
At the leaf, an entry of the matrix is emitted
for every row and column of the elemental matrix,
using the global row and column indices instead of the elemental ones.
No bottom-up phase is required for assembly,
as PETSc handles the merging of multi-instanced entries.






\section{Results} \label{sec: Results}

\subsubsection*{Computing Environment:\;}
We performed experiments, including simulations and scaling studies, on the Cascade Lake Compute Nodes of the \Frontera~ system. (Refer ~\secref{sec:Frontera} for compute configuration.)
\added{
\subsubsection*{Software and Libraries:\;} We used ~\petsc~ ~\cite{petsc-web-page} as the numerical algebra solver for solving system of equations. The \dsort~ and \tsort~ implementation is taken from \Dendro ~\cite{Dendro5}. All comparison with the immersed (IBM) method is performed using the open-source code~\cite{Dendrite} based on ~\citet{saurabh2020industrial}. Additionally, Matlab ~\cite{MATLAB:2018} is used for analyzing the condition number of matrices, and \textsc{trimesh} ~\cite{trimesh} is used to compute the signed distance. The roofline plot was generated by using Intel Advisor.
}

\subsection{Approximation of Voxelized Geometry}\label{sec: signeddistance}

\begin{figure}[b!]
\begin{subfigure}[b]{0.29\linewidth}
\includegraphics[trim=150 250 150 460, clip,width=\linewidth]{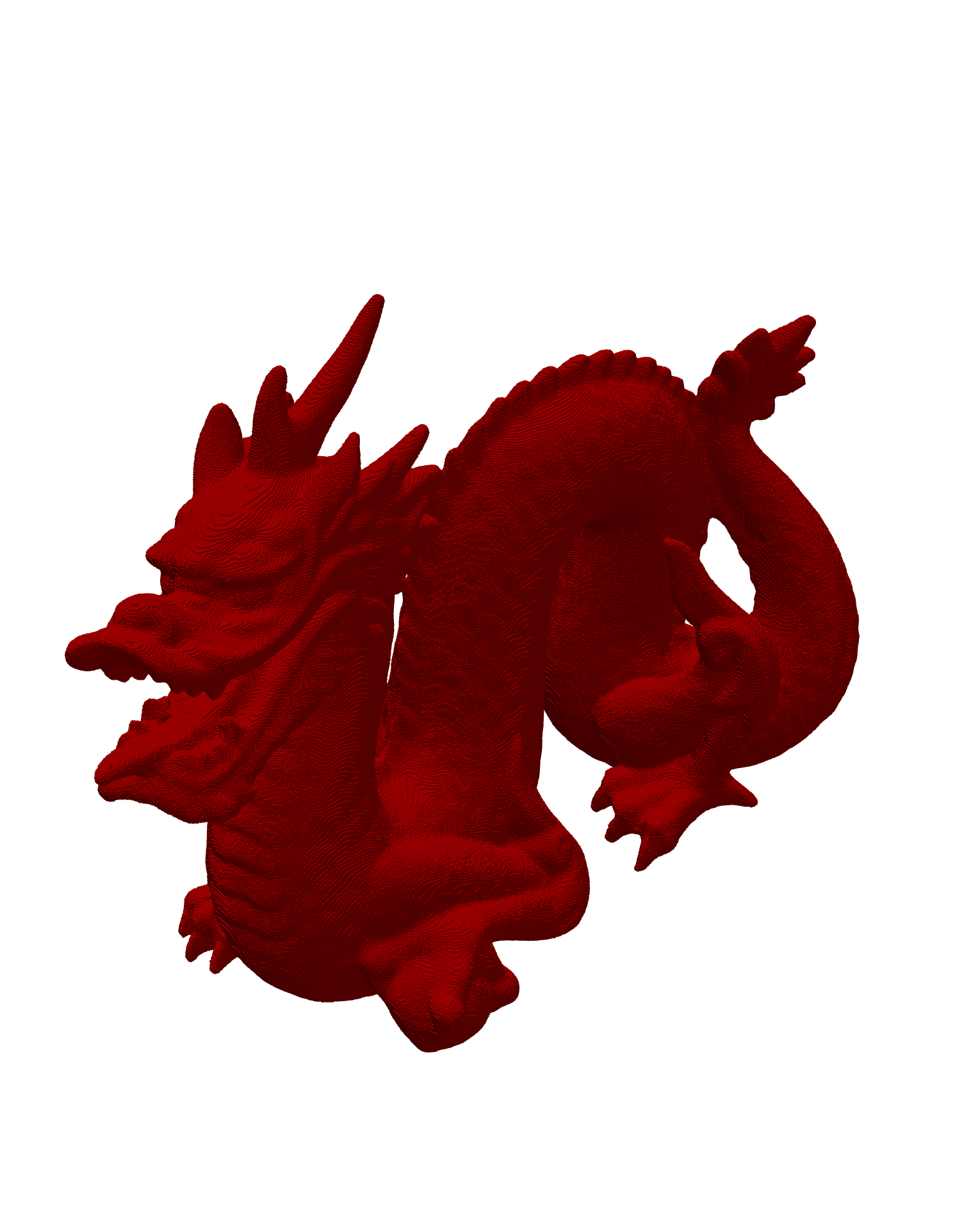}
\subcaption{Voxelized}
\label{fig:DragonVoxelized}
\end{subfigure}
\begin{subfigure}[b]{0.69\linewidth}
\centering
\begin{tikzpicture}[scale = 0.8]
  \begin{semilogyaxis}[
    width=1.0\linewidth, 
    height=0.65\linewidth, 
    axis y line*=left,
    xlabel= Refinement Level,
    ylabel= Number of elements,
  ]
  \addplot[mark=o,color = red,ultra thick] 
        table[x expr={\thisrow{RefineLvl}))},y expr=0.5*\thisrow{NumElements}]{Data/signedDistance.txt};
        \label{plot_2_y1}
    \end{semilogyaxis}

    \begin{semilogyaxis}[
      width=1.0\linewidth, 
      height=0.65\linewidth, 
      axis y line*=right,
      axis x line=none,
      ylabel=Distance,
      legend style={at={(0.5,1.4)},anchor=north,legend columns=2}, 
    ]
   \addplot[mark=o,color = blue,ultra thick] 
        table[x expr={\thisrow{RefineLvl}))},y expr={\thisrow{Distance}}]{Data/signedDistance.txt};\label{plot_2_y2} 
    \addlegendimage{/pgfplots/refstyle=plot_2_y1}\addlegendentry{Distance}
    \addlegendimage{/pgfplots/refstyle=plot_2_y2}\addlegendentry{Element count}
  \end{semilogyaxis}
\end{tikzpicture}
\subcaption{Error} 
\label{fig:DragonError}
\end{subfigure}
\vspace{0.15in}
\caption{Figure showing the voxelized geometry for the Stanford Dragon on octree mesh \figref{fig:DragonVoxelized}. The error (\figref{fig:DragonError}) is measured as the maximum of signed distance from boundary nodes of octree to the STL mesh. With increase in the refinement at the surface of geometry, the octree mesh coincides with the actual 3D mesh resulting in decrease in the signed distance error. Note the first order convergence in signed distance error with resolution.}
\label{fig: SignedDistance}
\end{figure}
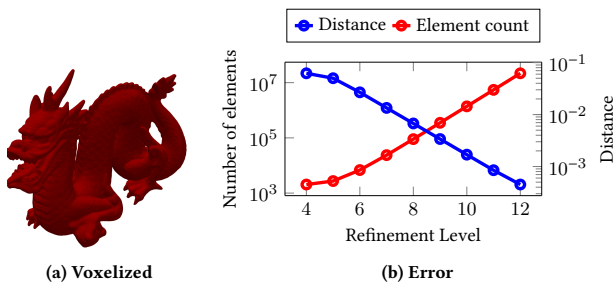

The \textit{carving-out} approach leads to a voxelized geometry, which is an approximation of the actual geometry. In this section, we compare how closely the voxelized geometry mimics the actual geometry by considering the example of the Stanford Dragon~\citep{levoy2005stanford}. 
\figref{fig: SignedDistance} compares the difference in the representation of actual boundary for the voxelized geometry by computing the signed distance\footnote{computed using \texttt{trimesh} library. A positive value denotes inside.}.  \figref{fig:DragonVoxelized} shows the voxel representation for the Stanford dragon. \figref{fig:DragonError} compares the $L_{\infty}$ error of computed signed distance between the boundary nodes of the voxelized geometry and the actual STL file. Similar to the previous case, we can see that with increase in the refinement, the voxelized geometry approaches the actual geometry.

\subsection{Conditioning of Discrete Operators} \label{sec:ConditionNumber}
As stated earlier, one approach to deploy traditional octrees on elongated channels is to stretch the mesh along the elongated channel ~\cite{esmaily2018scalable,mani2012analysis}. But this has a detrimental effect on condition number, which in turn will deteriorate the convergence of linear solvers. Table~\ref{table: conditionNumber} compares the variation in the condition number \footnote{evaluated with Matlab \texttt{condest} command} with the stretching of the elements for a  Laplace operator in 2D. We can see that with the increase in the aspect ratio of the mesh, the condition number of the linear system increases. With the generation of incomplete octree, we can ensure the aspect ratio of  each element in the mesh remains 1. Furthermore, since the error is dominated by the coarsest resolution, the incomplete octree permits decreasing the overall DOFs, at a given coarse resolution. This, in turn, decreases the condition number of the linear system.
\begin{table}[t!]
\newcommand{\tabincell}[2]{\begin{tabular}{@{}#1@{}}#2\end{tabular}}
\footnotesize
\begin{tabular}{|c|c|r|c|r|}
\hline
\multirow{2}{*}{\begin{tabular}[c]{@{}c@{}}Channel \\ length\end{tabular}} & \multicolumn{2}{c|}{Complete octree} & \multicolumn{2}{c|}{Incomplete octree} \\ \cline{2-5} 
 & DOFs & \tabincell{c}{Condition\\Number} & DOFs & \tabincell{c}{Condition\\Number} \\ \hline
1 & 1089 & 402.6 & 1089 & 402.6 \\ \hline
2 & 1089 & 466.7 & 561 & 155.6 \\ \hline
4 & 1089 & 510.1 & 297 & 42.5 \\ \hline
8 & 1089 & 512.0 & 165 & 13.3 \\ \hline
16 & 1089 & 10580.5 & 99 & 5.0 \\ \hline
\end{tabular}
\vspace{5 mm}
\caption{\footnotesize{Comparison of condition number for the case with complete octree and incomplete octree. In the case of complete octree, each element of the mesh was stretched according to the channel aspect ratio (represented here by the length) to conform with the channel boundaries, whereas in the case of incomplete octree, the aspect ratio was fixed to be 1 and the elements outside the domain are removed.}}
\label{table: conditionNumber}
\end{table}

\subsection{Convergence Test for Discrete Operators}
 Here, we present the convergence analysis for the Poisson operator $- \Delta u = f$ over the domain $\Omega$ with $u  = u_D$ on the domain boundary $\Gamma$. Inserting appropriate finite dimensional function spaces for trial and test function, the weak form of the Poisson operator can be written as: $(\nabla w^h,\nabla u^h)_{\Omega} = (w^h,f)_{\Omega}$
\begin{equation}
    (\nabla w^h,\nabla u^h)_{\Omega} = (w^h,f)_{\Omega}
\end{equation}

As mentioned previously, deploying incomplete octree based methods results in a voxelated geometry for a complicated geometrical shape. As discussed in~\secref{sec: signeddistance}, with the increase in the refinement of the element, the voxelated geometry approaches the true geometry. The rate of convergence of the distance follows only first order. Therefore, careful treatment is needed at the boundary elements to ensure an accurate order of convergence. In this context, several methods have been proposed in the literature~\cite{mittal2005immersed,lee2013implicit,burman2012fictitious}. In this work, we use the Shifted Boundary Method (SBM)~\cite{main2018shifted,atallah2020second} to treat boundary elements. 

The main idea behind SBM is to reformulate the original boundary value problem over a surrogate computational domain by modifying the original boundary conditions using Taylor series expansions. The weak form of Poisson operator after applying SBM treatment can be written as:
\begin{equation*}
\begin{split}
    (\nabla w^h, \nabla u^h)_{\Tilde{\Omega}} - (w^h,\nabla u^h \cdot \Tilde{\mathbf{n}})_{\Tilde{\Gamma}} - (\nabla w^h \cdot \Tilde{\mathbf{n}}, u^h + \nabla u^h \cdot \mathbf{d} -u_D)_{\Tilde{\Gamma}} + \\
    \frac{\alpha}{h}(w^h + \nabla w^h \cdot \mathbf{d}, u^h + \nabla u^h \cdot \mathbf{d} - u_D)_{\Tilde{\Gamma}} = (w^h,f)_{\Tilde{\Omega}}
\end{split}
\end{equation*}
where $\Tilde{\Omega}$ is the voxelated domain, $\Tilde{\Gamma}$ is the surface of the voxelated domain, $\Tilde{n}$ is the unit normal of the voxelated surface, $\alpha$ is the penalty term, $h$ is the element length, and $\textbf{d}$ is the distance vector from the boundary surface of voxelated domain $\Tilde{\Gamma}$ to the true surface $\Gamma$. The main idea of the method is to shift the boundary condition from $\Gamma$ to $\Tilde{\Gamma}$ by using second-order accurate Taylor series expansion. We omit the details here and refer to~\cite{main2018shifted,atallah2020second} for detailed analysis.

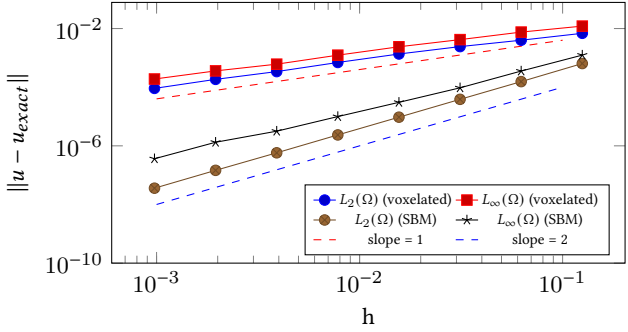
\begin{figure}[t!]
		\centering
		\begin{tikzpicture}[scale = 0.99]
		\begin{loglogaxis}[width=\linewidth, height=0.6\linewidth, scaled y ticks=true,xlabel={h},ylabel={$\norm{u - u_{exact}}$},
		legend entries={$L_2(\Omega)$ (voxelated),$L_\infty(\Omega)$ (voxelated),$L_2(\Omega)$ (SBM),$L_\infty(\Omega)$ (SBM),  slope = 1, slope = 2},
		legend style={at={(0.5,-0.25)},anchor=north, nodes={scale=0.65, transform shape}}, 
		legend pos= south east, 
		legend columns=2,
		ymin=1e-10,
		]
		\addplot table [x expr={1/2^\thisrow{h}},y={L2Omega},col sep=space] {Data/VoxelConvergence.txt};
		\addplot table [x expr={1/2^\thisrow{h}},y={LinfOmega},col sep=space] {Data/VoxelConvergence.txt};
		\addplot table [x expr={1/2^\thisrow{h}},y={L2Omega},col sep=space] {Data/convergence.txt};
		\addplot table [x expr={1/2^\thisrow{h}},y={LinfOmega},col sep=space] {Data/convergence.txt};
		\addplot +[mark=none, red, dashed] [domain=0.001:0.1]{0.04*x^1};
		\addplot +[mark=none, blue, dashed] [domain=0.001:0.1]{0.010*x^2};
		\end{loglogaxis}
		\end{tikzpicture}
		\caption{\textit{Convergence plot}: Figure showing the convergence behaviour for the Poisson operator on a two-dimensional circular disk. \vspace{0 mm}} 
		\label{fig:convergence}
	\end{figure}

To perform the convergence study, we consider Poisson problem on a two-dimensional disk of radius $R = 0.5$, centered at $(x_0 = 0.5,y_0 = 0.5)$ and $f = 1$. An exact solution exists and can be written as:
\begin{equation}
    u(r) = 0.25(R^2 - r^2)
\end{equation}
where $r = \sqrt{(x-x_0)^2+(y - y_0)^2}$.~\figref{fig:convergence} shows the convergence behaviour for the linear basis function. If we naively apply the boundary condition at the boundary nodes of the voxelated geometry, we only get a first-order convergence in both $L_2$ and $L_\infty$ norm. This is because the right boundary condition is applied at the wrong place, which is shifted by a distance $\mathbf{d}$ from the true boundary. As seen from the signed distance plot (\figref{fig:DragonError}), the voxelated geometry boundary approaches the true geometry according to the first order, and so is the convergence for the discrete Poisson operator. With the SBM method, we recover back the theoretical second-order convergence in both $L_2$ and $L_\infty$ norm for the linear basis function.

\subsection{Comparison with \textit{Immersed} Case}
Here, we present the comparison of the carved out approach with the immersed approach in terms of the  number of DOF and the total number of elements. We note that this analysis is equation agnostic. In order to compare the overall mesh element size and DOF, we set the background mesh to a constant refinement level and refined it near the object. \tabref{tab:dofCompar} compares the fraction of elements and DOF required for immersed and carved out approach. In the carved out case, all the elements and nodes that are inside the domain are discarded during the tree construction as mentioned in \secref{sec: octreeConstruct}, whereas for the immersed case, we retain the complete octree mesh. 2:1 balancing of octrees leads to the ripple effect, because of which there is a significant number of elements that are inside the domain (\figref{fig: meshComparison}). The nodes and elements that are marked \In (i.e. inside the object, sphere/dragon in this case) do not contribute towards the accuracy of the solution. These are not solved for in the system of equations, and eventually, a Dirichlet boundary condition is applied to it, but they had the associated cost during tree traversal and memory footprint for matrix and vector storage. \added{Overall, we see about an increase of 80--90\% in element size and a 33--40\% in the DOF count if we immerse an object.} The excess DOF count is significantly smaller than the element count because of the fact that we are performing continuous Galerkin (CG) computations and several elements share a common DOF. Additionally, we must recall the fact that in CG computations, hanging nodes do not contribute to the additional degrees of freedom. However, if we were to perform discontinuous Galerkin (DG) computation, each element would have its own unique node id and associated DOF. In such computations, the excess DOF count would scale as excess element count. The actual fraction of DOF and element that is reduced as a result of carving out depends upon the surface area and volume of the object that is being carved out. \added{A large surface area of the object would result in more elements at the finest resolution near the boundaries of the object. In contrast, a larger volume would result in more elements being discarded out from the interior of the object.} Constructing an incomplete octree by cutting the elements inside the object results in processing fewer elements during a solve. 
\added{
{\renewcommand{\arraystretch}{1.1} 
\begin{table}[t!]
{%
\footnotesize
{
\begin{tabular}{|c|c|c|c|c|c|}
\hline
\multicolumn{2}{|c|}{\multirow{2}{*}{}} & \multicolumn{4}{c|}{Refine Level} \\ \cline{3-6} 
\multicolumn{2}{|c|}{}  & 11 & 12 & 13 & 14 \\ \hline
\multirow{2}{*}{Sphere} & $f_{\text{elem}}$ & 1.75 & 1.79 & 1.81 & 1.82 \\ \cline{2-6} 
 & $f_\text{DOF}$ & 1.30 & 1.31 & 1.32 & 1.33 \\ \hline
\multirow{2}{*}{\begin{tabular}[c]{@{}c@{}}Stanford\\ Dragon\end{tabular}} & $f_{\text{elem}}$ & 1.84 & 1.87 & 1.90 & 1.92 \\ \cline{2-6} 
 & $f_\text{DOF}$ & 1.36 & 1.39  & 1.41 & 1.43 \\ \hline
\end{tabular}
}
}
\vspace{5mm}
\caption{\added{\footnotesize{Comparison of the ratio of number of elements ($f_\text{elem}$) and degrees of freedom ($f_\text{DOF}$) with and without (\textit{immersed}) carving out the sphere and the Stanford dragon from the domain. The base refinement was set to 4 and the refinement level near the object was varied from 11 to 14.}}}
\label{tab:dofCompar}
\end{table}
}
}

\subsection{Scaling}

\added{We evaluated the strong and weak scaling performance of our traversal-based \mvec~ using linear and quadratic elements on the \Frontera~ supercomputer for two different cases: a) an elongated channel of dimension $16 \times 1 \times 1$, b) a spherical region carved out from the cube. We individually timed the execution of major components of \mvec, namely top-down and bottom-up traversal, leaf \mvec~ to compute the elemental operators, malloc and communication cost. It must be noted that for any PDE solver, \mvec~ is the basic building block and determines the overall parallel performance and scalability. 
}
\added{We highlight some important points regarding the experimental setup for performing the scaling studies and the interpretability of the scaling results:
\begin{itemize}[left=0in]
    \item \textit{Strong Scaling:\;} For each of the strong scaling cases, we generated a fixed mesh defined by different refinement levels in different regions of interest. When comparing the mesh with linear and quadratic elements, the total number and distribution of elements in a given mesh is the same not only at a global level but also locally at each processor level. Note that the partitioning algorithm \dsort~ is agnostic to the underlying element order and distributes the element at the octant level  before the nodal information is encoded~\footnote{More formally, consider a mesh $M$ with $N$ global elements distributed over $p$ processor. If we globally number the elements of mesh from $0 \cdots N-1$, then if a processor $k, k\leq p$, receives $m_i \cdots m_k$ ($m_i$'s being the global element number)  elements for linear, $0 \leq m_i \leq m_k \leq N-1$, then the processor $k$ for quadratic mesh will also receive the same sequence of elements $m_i \cdots m_k$.}.
    \item However, the total number of DOF and problem size grows as $\mathcal{O}((p+1)^d)$ for an arbitrary order $p$ and dimension $d$~\footnote{Every element has $\mathcal{O}((p+1)^d)$ nodes (Refer \secref{createNodes}).}. Hence, the mesh with linear and quadratic elements have different computation and communication complexity. For instance, both linear and quadratic mesh for channel strong scaling study have 13.5M elements. But  the linear element mesh has 13.7M DOFs, whereas the quadratic has 109.1M DOFs.
    \item \textit{Weak Scaling:\;} For the weak scaling runs, with an increase in the number of processors, we increase the refinement level in the regions of interest in such a way that the average number of elements per processor remains the same. Similar to the strong scaling, for a given number of processors, the total number and distribution of elements are the same for both linear and quadratic basis functions both globally and locally. Hence, the quadratic mesh has a greater number of DOF compared to the linear one.
    \item In all the scaling figures (\figref{fig:strong-scaling} --\figref{fig:weak-scalingSphere}), for a given number of processors, the left bar corresponds to the \mvec~ execution profile for the linear elements, and the right bar corresponds to the execution profile of the quadratic elements. The total execution time for the linear elements is shown by solid blue lines and red dashed lines for the quadratic.
\end{itemize}
}
\subsubsection{Scaling results for the channel:\;}\label{sec:channelScaling}
The incomplete octrees representing $16 \times 1 \times 1$ elongated channel, with greater refinement on the boundary and minimal refinement on the interior, are generated to carry out the scaling studies. This is representative of the common cases that arise in the boundary-dominated physical phenomena. Each scaling run was repeated for linear and quadratic hexahedral grids.

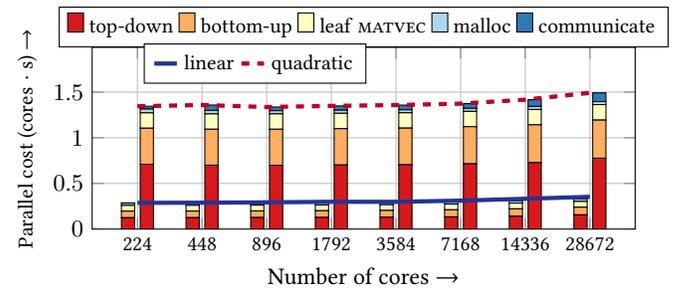
\begin{figure}[tbh]
  \centering
  \begin{tikzpicture}
  \begin{axis}[
  ybar stacked, bar width=5pt,    
  xlabel={Number of cores $\rightarrow$},
  ylabel={\small{Parallel cost  (cores $\cdot$ s) $\rightarrow$}},
   symbolic x coords={224,448,896,1792,3584,7168,14336,28672},width=8.8cm,height=4cm,
  xtick = data, 
  ymin=0,
  ymax = 1.99,
  x tick label style={font=\small},
  legend pos=north west,grid=major,legend style={at={(0.5,1.0)},
		anchor=south,legend columns=5}]

  \addplot [fill=div_d1] [bar shift=-.125cm] table[x={npes} , y expr={\thisrow{topdown(mean)}*\thisrow{npes}/2240} ]{Data/Channel/LinearSS.txt};
  \addplot [fill=div_d2] [bar shift=-.125cm]  table[x={npes}, y expr={\thisrow{bottomup(mean)}*\thisrow{npes}/2240} ]{Data/Channel/LinearSS.txt};
  \addplot [fill=div_d3] [bar shift=-.125cm]  table[x={npes}, y expr={\thisrow{elemental(mean)}*\thisrow{npes}/2240} ]{Data/Channel/LinearSS.txt};
  \addplot [fill=div_d4] [bar shift=-.125cm]  table[x={npes}, y expr={(\thisrow{matvec(mean)}-\thisrow{elemental(mean)}-\thisrow{bottomup(mean)}-\thisrow{topdown(mean)})*\thisrow{npes}/2240} ]{Data/Channel/LinearSS.txt};
  \addplot [fill=div_d5] [bar shift=-.125cm] table[x={npes}, y expr={\thisrow{ghostexchange(mean)}*\thisrow{npes}/2240} ]{Data/Channel/LinearSS.txt};
  
  \makeatletter
  \newcommand\resetstackedplotsA{
  \makeatletter
  \pgfplots@stacked@isfirstplottrue
  \makeatother
  \addplot [forget plot,draw=none] coordinates{(224,0) (448,0) (896,0) (1792,0) (3584,0) (7168,0) (14336,0) (28672,0)};
  }
  \makeatother
  \resetstackedplotsA
 \addplot [fill=div_d1] [bar shift=.125cm] table[x={npes} , y expr={\thisrow{topdown(mean)}*\thisrow{npes}/2240} ]{Data/Channel/QuadraticSS.txt};
  \addplot [fill=div_d2] [bar shift=.125cm]  table[x={npes}, y expr={\thisrow{bottomup(mean)}*\thisrow{npes}/2240} ]{Data/Channel/QuadraticSS.txt};
  \addplot [fill=div_d3] [bar shift=.125cm]  table[x={npes}, y expr={\thisrow{elemental(mean)}*\thisrow{npes}/2240} ]{Data/Channel/QuadraticSS.txt};
  \addplot [fill=div_d4] [bar shift=.125cm]  table[x={npes}, y expr={(\thisrow{matvec(mean)}-\thisrow{elemental(mean)}-\thisrow{bottomup(mean)}-\thisrow{topdown(mean)})*\thisrow{npes}/2240} ]{Data/Channel/QuadraticSS.txt};
  \addplot [fill=div_d5] [bar shift=.125cm] table[x={npes}, y expr={\thisrow{ghostexchange(mean)}*\thisrow{npes}/2240} ]{Data/Channel/QuadraticSS.txt};
  
  \legend{\small{top-down},\small{bottom-up}, \small{leaf \mvec}, \small{malloc}, \small{communicate}}
  \end{axis}
     \begin{axis}[ axis y line=none,axis x line = none,symbolic x coords={224,448,896,1792,3584,7168,14336,28672},width=8.8cm,height=4cm,
        xtick = data,ymin=0,ymax=1.99,legend style={at={(0.3,0.8)},
		anchor=south,legend columns=2}] 
      
      \addplot[sq_b1,ultra thick]  table[x={npes}, y expr={(\thisrow{matvec(mean)} + \thisrow{ghostexchange(mean)})*\thisrow{npes}/2240} ]{Data/Channel/LinearSS.txt};
      \addplot[sq_r1,ultra thick,dashed]  table[x={npes}, y expr={(\thisrow{matvec(mean)} + \thisrow{ghostexchange(mean)})*\thisrow{npes}/2240} ]{Data/Channel/QuadraticSS.txt};
      
      \legend{\small{linear},\small{quadratic}}
  
  \end{axis}
  \end{tikzpicture}
  \caption{\added{\textit{Strong scaling for channel case.\;} Parallel cost evaluated with the 3D Poisson \mvec~ on \Frontera~ supercomputer. Problem size was fixed at 13M elements (13.7M unknowns for linear and 109.1M unknowns for quadratic)}}
  \label{fig:strong-scaling}
\end{figure}

\added{For the strong scaling runs, we generated octree mesh with 13M elements for linear and quadratic basis functions. Both linear and quadratic mesh is similar at the elemental level. \figref{fig:strong-scaling}  shows the strong scaling behavior in terms of parallel cost (Run time $\times$ number of cores) for both the linear and quadratic basis functions. A constant line would mean ideal strong scaling efficiency. For the linear mesh, \mvec~ execution time decreased from 2.87 s on 224 processors to 0.027 s on 28K processors, resulting in 81\% parallel efficiency for 128 fold increase in processor count. Similarly, for the quadratic mesh, we see a reduction in \mvec~ execution time from 13.5 s on 224 processors to 0.1 s on 28K processors, resulting in 90\% parallel efficiency. The overall theoretical complexity for \mvec~ for a given element of order $p$ has been shown to scale as $\mathcal{O}( d(p+1)^{d+1})$. We see a factor of 4.2 $\times$ increase in \mvec~ execution time for quadratic element $(p=2)$ over linear $(p=1)$, which is within the theoretical bounds.}

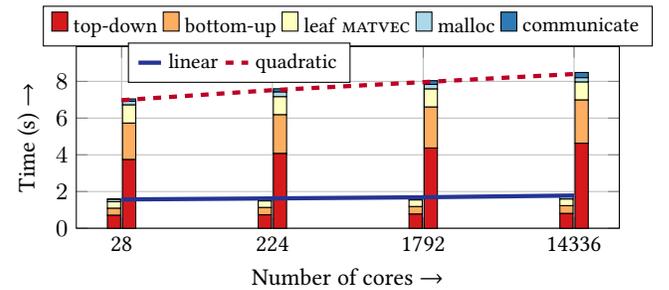
\begin{figure}[tbh]
  \centering
  \begin{tikzpicture}
  \begin{axis}[
  ybar stacked, bar width=5pt,    
  xlabel={Number of cores $\rightarrow$},
  ylabel={Time (s) $\rightarrow$ },
  ymin=0.0,
  ymax=9.9,
  symbolic x coords={28,224,1792,14336},width=8.8cm,height=4cm,
  xtick = data, 
  legend pos=north west,grid=major,legend style={at={(0.5,1.0)},
		anchor=south,legend columns=5}]

  \addplot [fill=div_d1] [bar shift=-.101cm] table[x={npes} , y expr={\thisrow{topdown(mean)}} ]{Data/Channel/LinearWS.txt};
  \addplot [fill=div_d2] [bar shift=-.101cm]  table[x={npes}, y expr={\thisrow{bottomup(mean)}} ]{Data/Channel/LinearWS.txt};
  \addplot [fill=div_d3] [bar shift=-.101cm]  table[x={npes}, y expr={\thisrow{elemental(mean)}} ]{Data/Channel/LinearWS.txt};
  \addplot [fill=div_d4] [bar shift=-.101cm]  table[x={npes}, y expr={(\thisrow{matvec(mean)}-\thisrow{elemental(mean)}-\thisrow{bottomup(mean)}-\thisrow{topdown(mean)})} ]{Data/Channel/LinearWS.txt};
  \addplot [fill=div_d5] [bar shift=-.101cm] table[x={npes}, y expr={\thisrow{ghostexchange(mean)}} ]{Data/Channel/LinearWS.txt};
  
  \makeatletter
  \newcommand\resetstackedplotsThree{
  \makeatletter
  \pgfplots@stacked@isfirstplottrue
  \makeatother
   \addplot [forget plot,draw=none] coordinates{(28,0) (224,0) (1792,0) (14336,0)};
 }
  \makeatother
  \resetstackedplotsThree
  \addplot [fill=div_d1] [bar shift=+.101cm] table[x={npes} , y expr={\thisrow{topdown(mean)}} ]{Data/Channel/QuadWS.txt};
  \addplot [fill=div_d2] [bar shift=+.101cm]  table[x={npes}, y expr={\thisrow{bottomup(mean)}} ]{Data/Channel/QuadWS.txt};
  \addplot [fill=div_d3] [bar shift=+.101cm]  table[x={npes}, y expr={\thisrow{elemental(mean)}} ]{Data/Channel/QuadWS.txt};
  \addplot [fill=div_d4] [bar shift=+.101cm]  table[x={npes}, y expr={(\thisrow{matvec(mean)}-\thisrow{elemental(mean)}-\thisrow{bottomup(mean)}-\thisrow{topdown(mean)})} ]{Data/Channel/QuadWS.txt};
  \addplot [fill=div_d5] [bar shift=+.101cm] table[x={npes}, y expr={\thisrow{ghostexchange(mean)}} ]{Data/Channel/QuadWS.txt};
  \legend{\small{top-down}, \small {bottom-up}, \small {leaf \mvec},\small {malloc}, \small {communicate}}
  \end{axis}
  
    \begin{axis}[ axis y line=none,axis x line = none,symbolic x coords={28,224,1792,14336},width=8.8cm,height=4cm,
        xtick = data,ymin=0,ymax=10.0,legend style={at={(0.3,0.8)},
		anchor=south,legend columns=2}] 
      
      \addplot[sq_b1,ultra thick]  table[x={npes}, y expr={(\thisrow{matvec(mean)} + \thisrow{ghostexchange(mean)})} ]{Data/Channel/LinearWS.txt};
      \addplot[sq_r1,ultra thick,dashed]  table[x={npes}, y expr={(\thisrow{matvec(mean)} + \thisrow{ghostexchange(mean)})} ]{Data/Channel/QuadWS.txt};
      \legend{\small{linear},\small{quadratic}}
  \end{axis}  
  \end{tikzpicture}
  
  \caption{\textit{Weak scaling run time for channel case:\;} Execution time of 3D Poisson \mvec~ on \Frontera~ supercomputer, for a fixed grain size of about 35K elements per core. }
  \label{fig:weak-scaling}
\end{figure}

\added{
For the weak scaling runs, we created grids with a fixed grain size of about 35K elements per core and timed \mvec~ execution time. The coarsest mesh consists of 981K elements on 28 processors  with 1.02M DOFs for linear and 8.01M DOFs for quadratic element, whereas the finest mesh consists of 502M elements on a 14K processors with 505M DOFs for linear and 4 billion DOFs for quadratic elements.  \figref{fig:weak-scaling} plots the mean execution of the  \mvec~ averaged over 100 iterations as a function of the number of cores. A constant execution time would imply ideal weak scaling efficiency. We observed a slowly growing weak-scaled execution time. Overall the time increased from about 1.58 s on 28 cores to 1.9 s on 14 K cores for linear elements (82\% weak scaling efficiency) and 7.04 s to 8.04 s for quadratic elements (86\% weak scaling efficiency).}

\subsubsection{Scaling results for a spherical carved out region:\; }\label{sec:sphereScaling} To study the scaling behavior for a complex carved out geometry, we carved out a spherical region from a cubical domain. A sphere of diameter $d = 1$ unit is carved out from a cubical domain of $10 \times 10 \times 10$. Overall, full mesh contains 5 levels of octree adaptivity with maximum refinement near the sphere. Such domain and mesh resolution are similar to the application problem used for validation of Navier--Stokes simulation. 

\begin{figure}[tbh]
  \centering
  \begin{tikzpicture}
  \begin{axis}[
  ybar stacked, bar width=5pt,    
  xlabel={Number of cores $\rightarrow$},
  ylabel={\small{Parallel cost (cores $\cdot$ s) $\rightarrow$} },symbolic x coords={224,448,896,1792,3584,7168},width=9cm,height=4cm,
  xtick = data, 
  ymin=0,
  ymax=2.8,
  legend pos=north west,grid=major,legend style={at={(0.5,1.0)},
		anchor=south,legend columns=5}]

  \addplot [fill=div_d1] [bar shift=-.125cm] table[x={npes} , y expr={\thisrow{topdown(mean)}*\thisrow{npes}/2240} ]{Data/SphereScaling/Strong/sphereLinear.txt};
  \addplot [fill=div_d2] [bar shift=-.125cm]  table[x={npes}, y expr={\thisrow{bottomup(mean)}*\thisrow{npes}/2240} ]{Data/SphereScaling/Strong/sphereLinear.txt};
  \addplot [fill=div_d3] [bar shift=-.125cm]  table[x={npes}, y expr={\thisrow{elemental(mean)}*\thisrow{npes}/2240} ]{Data/SphereScaling/Strong/sphereLinear.txt};
  \addplot [fill=div_d4] [bar shift=-.125cm]  table[x={npes}, y expr={((\thisrow{matvec(mean)}-\thisrow{elemental(mean)}-\thisrow{bottomup(mean)}-\thisrow{topdown(mean)})*\thisrow{npes})/2240} ]{Data/SphereScaling/Strong/sphereLinear.txt};
  \addplot [fill=div_d5] [bar shift=-.125cm] table[x={npes}, y expr={\thisrow{ghostexchange(mean)}*\thisrow{npes}/2240} ]{Data/SphereScaling/Strong/sphereLinear.txt};
  
  \makeatletter
  \newcommand\resetstackedplotsSix{
  \makeatletter
  \pgfplots@stacked@isfirstplottrue
  \makeatother
  \addplot [forget plot,draw=none] coordinates{(224,0) (448,0) (896,0) (1792,0) (3584,0) (7168,0)};
 }
  \makeatother
  \resetstackedplotsSix
  \addplot [fill=div_d1] [bar shift=.125cm] table[x={npes} , y expr={\thisrow{topdown(mean)}*\thisrow{npes}/2240} ]{Data/SphereScaling/Strong/sphereQuadratic.txt};
  \addplot [fill=div_d2] [bar shift=.125cm]  table[x={npes}, y expr={\thisrow{bottomup(mean)}*\thisrow{npes}/2240} ]{Data/SphereScaling/Strong/sphereQuadratic.txt};
  \addplot [fill=div_d3] [bar shift=.125cm]  table[x={npes}, y expr={\thisrow{elemental(mean)}*\thisrow{npes}/2240} ]{Data/SphereScaling/Strong/sphereQuadratic.txt};
  \addplot [fill=div_d4] [bar shift=.125cm]  table[x={npes}, y expr={((\thisrow{matvec(mean)}-\thisrow{elemental(mean)}-\thisrow{bottomup(mean)}-\thisrow{topdown(mean)})*\thisrow{npes})/2240} ]{Data/SphereScaling/Strong/sphereQuadratic.txt};
  \addplot [fill=div_d5] [bar shift=.125cm] table[x={npes}, y expr={\thisrow{ghostexchange(mean)}*\thisrow{npes}/2240} ]{Data/SphereScaling/Strong/sphereQuadratic.txt};
  
  \legend{\small{top-down}, \small{bottom-up}, \small{leaf \mvec}, \small{malloc}, \small{communicate}}
  \end{axis}
      
      
  
    \begin{axis}[ axis y line=none,axis x line = none,symbolic x coords={224,448,896,1792,3584,7168},width=9cm,height=4cm,xtick = data,legend pos={north west},legend columns=2,
    ymin=0,
  ymax=2.8] 
      \addplot[sq_b1,ultra thick]  table[x={npes}, y expr={((\thisrow{matvec(mean)} + \thisrow{ghostexchange(mean)})*\thisrow{npes}) /(2240)} ]{Data/SphereScaling/Strong/sphereLinear.txt};
      \addplot[sq_r1,ultra thick,dashed]  table[x={npes}, y expr={((\thisrow{matvec(mean)} + \thisrow{ghostexchange(mean)})*\thisrow{npes}) /(2240)} ]{Data/SphereScaling/Strong/sphereQuadratic.txt};
      \legend{\small{linear},\small{quadratic}}
  \end{axis}  
  \end{tikzpicture}
  \caption{\textit{Strong scaling for sphere case:\;} Parallel cost evaluated on \Frontera~ supercomputer. Problem size was fixed at 17.5 M elements(17.4M unknowns for linear and 139.7M unknowns for quadratic)}
  \label{fig:strong-scalingSphere}
\end{figure}
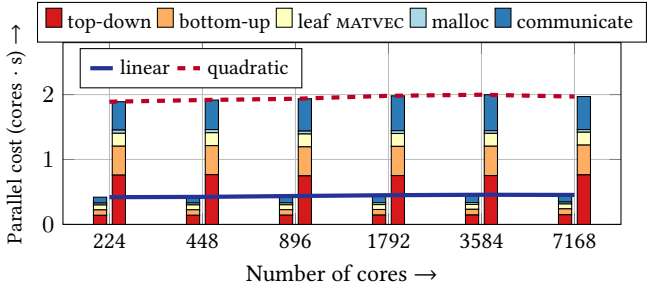

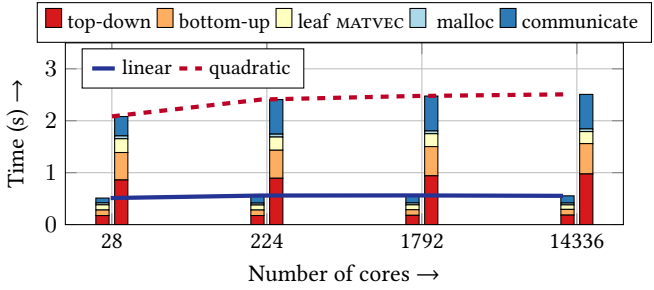
\begin{figure}[tbh]
  \centering
  \begin{tikzpicture}
  \begin{axis}[
  ybar stacked, bar width=5pt,    
  xlabel={Number of cores $\rightarrow$},
  ylabel={Time (s) $\rightarrow$ },symbolic x coords={28,224,1792,14336},width=9cm,height=4cm,
  xtick = data, 
  ymin=0,
  ymax=3.5,
  legend pos=north west,grid=major,legend style={at={(0.5,1.0)},
		anchor=south,legend columns=5}]

  \addplot [fill=div_d1] [bar shift=-.125cm] table[x={npes} , y expr={\thisrow{topdown(mean)}} ]{Data/SphereScaling/Weak/sphereLinear.txt};
  \addplot [fill=div_d2] [bar shift=-.125cm]  table[x={npes}, y expr={\thisrow{bottomup(mean)}} ]{Data/SphereScaling/Weak/sphereLinear.txt};
  \addplot [fill=div_d3] [bar shift=-.125cm]  table[x={npes}, y expr={\thisrow{elemental(mean)}} ]{Data/SphereScaling/Weak/sphereLinear.txt};
  \addplot [fill=div_d4] [bar shift=-.125cm]  table[x={npes}, y expr={((\thisrow{matvec(mean)}-\thisrow{elemental(mean)}-\thisrow{bottomup(mean)}-\thisrow{topdown(mean)}))} ]{Data/SphereScaling/Weak/sphereLinear.txt};
  \addplot [fill=div_d5] [bar shift=-.125cm] table[x={npes}, y expr={\thisrow{ghostexchange(mean)}} ]{Data/SphereScaling/Weak/sphereLinear.txt};
  
  \makeatletter
  \newcommand\resetstackedplotsSeven{
  \makeatletter
  \pgfplots@stacked@isfirstplottrue
  \makeatother
  \addplot [forget plot,draw=none] coordinates{(28,0) (224,0) (1792,0) (14336,0)};
 }
  \makeatother
  \resetstackedplotsSeven
 \addplot [fill=div_d1] [bar shift=.125cm] table[x={npes} , y expr={\thisrow{topdown(mean)}} ]{Data/SphereScaling/Weak/sphereQuadratic.txt};
  \addplot [fill=div_d2] [bar shift=.125cm]  table[x={npes}, y expr={\thisrow{bottomup(mean)}} ]{Data/SphereScaling/Weak/sphereQuadratic.txt};
  \addplot [fill=div_d3] [bar shift=.125cm]  table[x={npes}, y expr={\thisrow{elemental(mean)}} ]{Data/SphereScaling/Weak/sphereQuadratic.txt};
  \addplot [fill=div_d4] [bar shift=.125cm]  table[x={npes}, y expr={((\thisrow{matvec(mean)}-\thisrow{elemental(mean)}-\thisrow{bottomup(mean)}-\thisrow{topdown(mean)}))} ]{Data/SphereScaling/Weak/sphereQuadratic.txt};
  \addplot [fill=div_d5] [bar shift=.125cm] table[x={npes}, y expr={\thisrow{ghostexchange(mean)}} ]{Data/SphereScaling/Weak/sphereQuadratic.txt};
  
  \legend{\small{top-down},\small{bottom-up},\small{leaf \mvec},\small{ malloc},\small{communicate}}
  \end{axis}
    \begin{axis}[ axis y line=none,axis x line = none,symbolic x coords={28,224,1792,14336},width=8.8cm,height=4cm,xtick = data,legend pos={north west},legend columns=2,
    ymin=0,
  ymax=3.5] 
      \addplot[sq_b1,ultra thick]  table[x={npes}, y expr={((\thisrow{matvec(mean)} + \thisrow{ghostexchange(mean)})} ]{Data/SphereScaling/Weak/sphereLinear.txt};
      \addplot[sq_r1,ultra thick,dashed]  table[x={npes}, y expr={((\thisrow{matvec(mean)} + \thisrow{ghostexchange(mean)})} ]{Data/SphereScaling/Weak/sphereQuadratic.txt};
      \legend{\small{linear},\small{quadratic}}
  \end{axis}  
  \end{tikzpicture}
  \caption{\textit{Weak scaling run time for sphere case:\;} Mean Execution time of 100 \mvec~  on Frontera supercomputer, for a fixed grain size of about 10K elements per core.}
  \label{fig:weak-scalingSphere}
\end{figure}
{\renewcommand{\arraystretch}{1.1} 
\begin{table}[t!]
\resizebox{\linewidth}{!}{%
{
\Large
\begin{tabular}{|c|c|c|c|c|c|c|}
\hline
\multirow{2}{*}{\textbf{\begin{tabular}[c]{@{}c@{}}Case\\ Type\end{tabular}}} & \multirow{2}{*}{\textbf{\begin{tabular}[c]{@{}c@{}}Element\\ Order\end{tabular}}} & \multicolumn{3}{c|}{\textbf{Strong scaling}} & \multicolumn{2}{c|}{\textbf{Weak Scaling}} \\ \cline{3-7} 
 &  & \textbf{\begin{tabular}[c]{@{}c@{}}Num \\ elements\end{tabular}} & \textbf{\begin{tabular}[c]{@{}c@{}}Num\\ DOFs\end{tabular}} & \textbf{Efficiency} & \textbf{\begin{tabular}[c]{@{}c@{}}Num \\ elements/core\end{tabular}} & \textbf{Efficiency} \\ \hline
\multirow{2}{*}{\textbf{\begin{tabular}[c]{@{}c@{}}Channel\\ (\secref{sec:channelScaling})\end{tabular}}} & \textbf{Linear} & 13.5 M & 13.7 M & 0.81 & 35K & 0.82 \\ \cline{2-7} 
 & \textbf{Quadratic} & 13.5 M
 & 101.9 M & 0.90 & 35K & 0.86 \\ \hline
\multirow{2}{*}{\textbf{\begin{tabular}[c]{@{}c@{}}Sphere\\ (\secref{sec:sphereScaling})\end{tabular}}} & \textbf{Linear} & 17.5 M & 17.4 M & 0.90 & 10K & 0.74 \\ \cline{2-7} 
 & \textbf{Quadratic} & 17.5 M & 139.7 M & 0.96 & 10K & 0.83 \\ \hline
\end{tabular}
}
}
\vspace{5 mm}
\caption{\added{\footnotesize{Summary of scaling efficiency for the channel and spherical carved out region.}}}
\label{tab:scalCompar}
\end{table}
}
\added{We created grids of about 17.5M elements for the strong scaling, which correspond to 17.3M DOFs for the linear and 139.7M DOFs for the quadratic basis functions. \added{Similar to the channel case, the mesh partition for both the linear and quadratic is similar at the elemental level but has a different number of DOFs.} \figref{fig:strong-scalingSphere} shows the parallel efficiency averaged over 100 \mvec~ iterations. Overall we observe a good overall parallel efficiency. In the case of linear elements, we observe a 29$\times$ reduction in \mvec~ execution time for 32 fold increase in processor (90\% strong scaling efficiency). In contrast,  the quadratic element resulted in a 31$\times$ reduction in computation time (96\% strong scaling efficiency).}

For the weak scaling, we kept a constant grain size of around 10K elements per processor. The coarsest mesh consists of about 290K elements resulting in 280K DOFs for the linear basis function and 2.3 M for the quadratic. In contrast, the finest mesh consists of 138 M  elements with about 138M DOFs for linear and 1.1 billion DOFs for quadratic basis function. \figref{fig:weak-scalingSphere} shows the overall weak scaling performance. \mvec~ execution time grew from for 4.1 s on 28 processors to about 5.5 s on 14K processors for linear elements, resulting in a factor of about 1.34$\times$ increase for 512$\times$ in the number of processors (74\% efficiency). In the case of quadratic elements, the execution time increased from 20 s on 28 processors to about 25 s on 14K processors, yielding about 83\% weak scaling efficiency. \tabref{tab:scalCompar} summarizes the scaling efficiency for both the channel and the sphere case.
\begin{figure}[t!]
    \centering
    \begin{tikzpicture}[scale=0.7]
    \begin{loglogaxis}[ width=0.99\linewidth, height = 0.8\linewidth, axis y line*=left,legend pos=north east,legend style={at={(1.0,1.15)},legend columns=2},ylabel = Average number of ghost nodes,xlabel ={number of cores $\rightarrow$}]
        \addplot [blue,thick,error bars/.cd, y dir=both, y explicit]
            table [col sep=space, x=npes, y=Communicated(mean)    , y error=Communicated(std)]{Data/SphereScaling/Statistics/StatisticsLinear.txt};
            \addplot [red,thick,error bars/.cd, y dir=both, y explicit]
            table [col sep=space, x=npes, y=Communicated(mean)    , y error=Communicated(std)]{Data/SphereScaling/Statistics/StatisticsQuadratic.txt};
                \legend{\small {Linear}, \small{Quadratic (Number of ghost nodes)}}
    \end{loglogaxis}
    \begin{semilogxaxis}[ width=0.99\linewidth, height=0.8\linewidth, axis y line*=right,legend pos=north west,legend style={legend columns=1},axis x line=none,legend style={at={(0.28,0.95)},legend columns=2},
    ymin = 0,
    ymax = 1,
    ylabel ={$\eta = \frac{\mathrm{Ghost \; nodes}}{\mathrm{Owned \; nodes}}$}]
        \addplot [ densely dashed, ultra thick,error bars/.cd, y dir=both, y explicit,error bar style={solid}]
            table [col sep=space, x=npes, y=Ratio(mean)    , y error=Ratio(std)]
            {Data/SphereScaling/Statistics/StatisticsLinear.txt};
            \addplot [ densely dashed, brown,ultra thick,error bars/.cd, y dir=both, y explicit,error bar style={solid}]
            table [col sep=space, x=npes, y=Ratio(mean)    , y error=Ratio(std)]{Data/SphereScaling/Statistics/StatisticsQuadratic.txt};
                \legend{\small Linear, \small Quadratic ($\eta$)}
    \end{semilogxaxis}
    \end{tikzpicture}
    
    \caption{Figure showing the mean and standard deviation of the distribution of ghost nodes (shown by solid lines) and ratio of ghost nodes by owned nodes per processor (shown by dashed lines)\vspace{2 mm}}
    \label{fig: DAStatistics}
\end{figure}
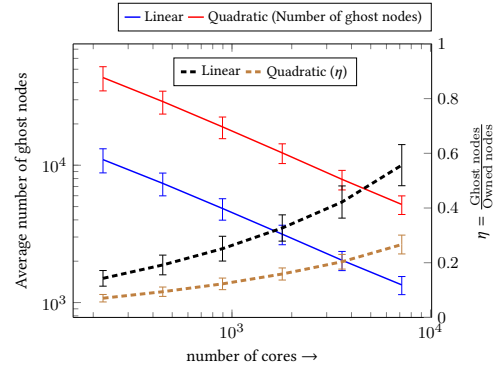


Further, we analyzed the distribution of ghost nodes per processor for the above sphere case, which is indicative of bytes of data communicated. In our experiment, we kept a similar distribution of the elements across processors for both linear and quadratic basis function but has different degrees of freedom associated with them. The amount of data communicated across processors is a function of the total number of ghost elements that share partition boundaries. With the increase in the number of processors, the total number of ghost elements increases, but the average number of ghost elements decreases. For an arbitrary order element $p$, the number of nodes that share faces across the processor boundaries (and hence needs exchange of information) grows as $\mathcal{O}((p+1)^{d - 1})$. Since the partition is similar at the elemental level, the average number of ghost nodes that are needed for ghost exchange  is higher for the quadratic compared to the linear elements. The solid lines in \figref{fig: DAStatistics} show the comparison for linear and quadratic case.

\added{Additionally, we also analyzed the distribution of the ratio of ghost nodes to the number of owned nodes (denoted by $\eta$), which is indicative of the extent over which the communication can be overlapped with computation. Let $N_L$ be the number of local nodes that are owned by processor (do not share processor boundaries with any other elements) and $N_G$ be the number of ghost nodes, then:
\begin{equation*}
    \eta = \frac{N_G}{N_L} \propto \frac{(p+1)^{d-1}}{(p+1)^{d}}  = \frac{1}{(p+1)}
\end{equation*}
From the above equation, we can see that this ratio grows inversely with respect to the degree of element. We observe similar behavior in our experiments, as shown by dashed lines in \figref{fig: DAStatistics}. This explains the better scaling efficiency for quadratic as compared to linear shown in \tabref{tab:scalCompar}. For a single processor run, $\eta = 0$. With an increase in the number of processors $\eta$ increases, and in extreme limit of parallelization, when each processor has only one element, $\eta \to 1$. It is non-trivial to analyze the exact rate of increase in $\eta$ for arbitrary shapes as a function of processor count and is beyond the scope of the current work.
}
\subsubsection{Roofline:\;}
\begin{figure}
\begin{tikzpicture}
\begin{loglogaxis}[
    log ticks with fixed point,
    width=0.99\linewidth,height=5cm,
    xmax = 2,
    xmin = 0.001,
    ymin = 0.1,
    xlabel = Arithmetic Intensity (AI),
    ylabel = GFLOP/s
]
\addplot [mark=square,mark options=solid,mark size = 1pt]table {
0.018 6.56
0.1 38.13
0.21 76.26
10 76.26
};
\addplot []table {
0.001 0.3685
0.018 6.56
};

\addplot [dotted,mark=square,mark options=solid,mark size = 1pt]table {
0.036 6.56
0.21 38.13
0.42 76.27
10 76.27
};
\addplot [dotted,mark=none]table {
0.001 0.18351
0.036 6.56
};

\addplot [dotted,mark=square,mark options=solid,mark size = 1pt]table {
0.26 6.56
1.5 38.13
3.01 76.27
};
\addplot [dotted,mark=none]table {
0.001 0.02534
0.26 6.56
};

\addplot [mark=square,mark options=solid,mark size = 1pt]table {
0.43 6.56
2.52 38.13
5.04 76.27
};
\addplot [mark=none]table {
0.001 0.01515
0.43 6.56
};

\addplot [mark=none]table {
0.018 6.56
0.036 6.56
0.26 6.56
0.43 6.56
10 6.56
};

\addplot [mark=none]table {
0.1 38.13
0.21 38.13
1.5 38.13
2.52 38.13
10 38.13
};

\draw[-] (axis cs:0.001,0.37)  edge  node[sloped, anchor=left, above, text width=2.2cm] {\fontsize{4pt}{4pt}\selectfont{\textbf{L1\;Bandwidth:\;368.5\;GB/s}}} (axis cs:0.018,6.56);
\draw[dotted] (axis cs:0.001,0.18351)  edge  node[sloped, above, text width=2.2cm] {\fontsize{4pt}{4pt}\selectfont{\textbf{L2\;Bandwidth:\;183.51\;GB/s}}} (axis cs:0.036,6.56);
\draw[dotted] (axis cs:0.001,0.02534)  edge  node[sloped, above, text width=2.2cm] {\fontsize{4pt}{4pt}\selectfont{\textbf{L3\;Bandwidth:\;25.34\;GB/s}}} (axis cs:0.26,6.56);
\draw[dotted] (axis cs:0.01,0.1515)  edge  node[sloped, below, text width=2.4cm] {\fontsize{4pt}{4pt}\selectfont{\textbf{DRAM\;Bandwidth:\;15.15\;GB/s}}} (axis cs:0.43,6.56);
\draw[dotted] (axis cs:0.2,76.26)  edge  node[sloped, above, text width=2.4cm] {\fontsize{4pt}{4pt}\selectfont{\textbf{Vector\;FMA\;Peak:\;76.26\;GFLOP/s}}} (axis cs:2.0,76.26);
\draw[dotted] (axis cs:0.2,38.13)  edge  node[sloped, above, text width=2.4cm] {\fontsize{4pt}{4pt}\selectfont{\textbf{Vector\;Add\;Peak:\;36.87\;GFLOP/s}}} (axis cs:2.0,38.13);
\draw[-] (axis cs:0.2,6.56)  edge  node[sloped, above, text width=2.4cm] {\fontsize{4pt}{4pt}\selectfont{\textbf{Scalar\;Add\;Peak:\;6.55\;GFLOP/s}}} (axis cs:2.0,6.56);
\draw[dashed,ao] (axis cs:0.001,0.060)  edge  node[sloped, above, text width=2.2cm] {\fontsize{4pt}{4pt}\selectfont{\textbf{Achieved\;Bandwidth:\;60GB/s}}} (axis cs:0.08,4.8);
\draw[dashed,ao] (axis cs:0.001,0.060) --  (axis cs:1.0,60.0);

\node [isosceles triangle,draw,red,thick,minimum size=1pt,scale = 0.3] at (axis cs:0.121,6.84) {};
\node [circle,draw,red,thick, minimum size=1pt,scale = 0.3] at (axis cs:0.121,5.557){};

\node [isosceles triangle,draw,blue,thick,minimum size=1pt,scale = 0.3] at (axis cs:0.0012,90) {}; 
\node [circle,draw,blue,thick, minimum size=1pt,scale = 0.3] at (axis cs:0.0012,60){};
\node [anchor=west] at (axis cs:0.00125,60) {\tiny{Linear: 65K elements}};
\node [anchor=west] at (axis cs:0.00125,90) {\tiny{Linear: 525K elements}};

\node [isosceles triangle,draw,red,thick,minimum size=1pt,scale = 0.3] at (axis cs:0.012,90) {}; 
\node [circle,draw,red,thick, minimum size=1pt,scale = 0.3] at (axis cs:0.012,60){};
\node [anchor=west] at (axis cs:0.0125,60) {\tiny{Quad: 65K elements}};
\node [anchor=west] at (axis cs:0.0125,90) {\tiny{Quad: 525K elements}};


\node [isosceles triangle,draw,blue,thick,minimum size=1pt,scale = 0.3] at (axis cs:0.072,4.261) {};
\node [circle,draw,blue,thick, minimum size=1pt,scale = 0.3] at (axis cs:0.072,3.646){};
\end{loglogaxis}
\end{tikzpicture}
\caption{\added{Figure showing roofline plot for the Poisson ~\mvec~ for linear and quadratic basis function for two different meshes on ~\Frontera. The plot was generated using Intel Advisor. The green dashed line shows the achieved bandwidth from our code. All values reported in the plot corresponds to double precision floating point operations.}}
\label{fig:roofline}
\end{figure}
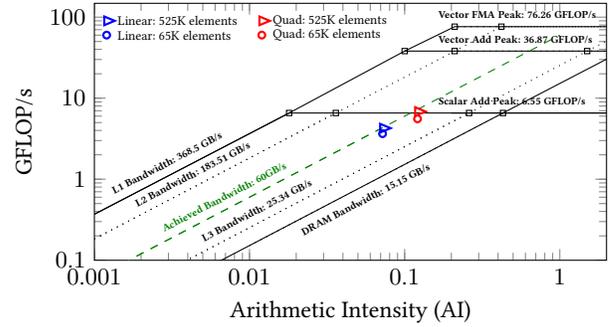
\added{~\figref{fig:roofline} shows the single core roofline plot for the elemental ~\mvec~ computation of Poisson operator using linear and quadratic basis function on ~\Frontera. Overall, we can see that the code is memory bound as is common for finite element codes. We observe higher arithmetic intensity \footnote{AI is measured as the amount of floating point operations performed per byte of data loaded into the memory.}
(AI) for quadratic (0.121) as compared to linear (0.072) elements. The amount of data needed for ~\mvec~ computation grows as $\mathcal{O}((p+1)^d)$ whereas the ~\mvec~ computation complexity grows as $\mathcal{O}(d(p+1)^{d+1})$. Therefore, AI tends to increase with polynomial order, which explains the observed behavior. We are able to achieve a performance of about 4 GFLOP/s using linear basis function and 7 GFLOP/s using quadratic basis functions for two different meshes, which corresponds to a bandwidth of approximately 60 GB/s as shown by the green lines. We note that we have not used any hand-coded explicit vectorization to ensure the portability of the code across various platforms and relied on compiler-directed vectorization. We would like to explore some future avenues from the code optimization point of view pertaining to more efficient cache blocking techniques and architecture-specific efficient vectorized implementation of tensor products.
}

\subsection{Comparison with Existing Method}\label{sec: comparD4}
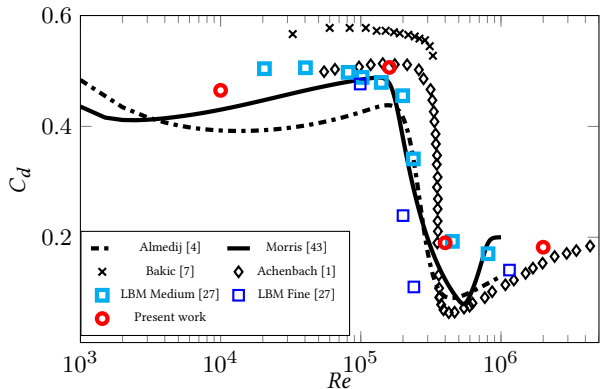
\begin{figure}
        
    \begin{subfigure}{1.0\linewidth}
    \hspace{0 mm}
    \centering
 \begin{tikzpicture}
      \begin{semilogxaxis}[
      ylabel near ticks,
      xlabel near ticks,
      ylabel shift=-2pt,
      xlabel shift=-6pt,
          width=1.0\linewidth, 
          height=0.7\textwidth,
          xlabel=$Re$, 
          ylabel=$C_d$,
         ymax = 0.6,
          ymin =  1E-2,
          xmin = 1000,
          xmax = 5E6,
        legend style={at={(0,0.02)},anchor= south west,legend columns=2}, 
        ]
        \addplot[mark=none,color = black,ultra thick,dashdotted] 
        table[x expr={\thisrow{Re}))},y expr=\thisrow{Drag},col sep=comma]{Data/DragCrisis/Morris.txt};
        \addplot[mark=none,color = black,ultra thick] 
        table[x expr={\thisrow{Re}))},y expr=\thisrow{Drag},col sep=comma]{Data/DragCrisis/Almedij.txt};
        \addplot[only marks,mark=x,mark size=2pt,color=black,thick] 
        table[x expr={\thisrow{Re}))},y expr=\thisrow{Drag},col sep=comma]{Data/DragCrisis/Bakic.txt};
        \addplot[only marks,mark=diamond,mark size=2pt,color=black,thick] 
        table[x expr={\thisrow{Re}))},y expr=\thisrow{Drag},col sep=comma]{Data/DragCrisis/Achenbach.txt};
        \addplot[only marks,mark=square ,mark size=2pt,color=cyan,ultra thick]
        table[x expr={\thisrow{Re}))},y expr=\thisrow{Drag},col sep=comma]{Data/DragCrisis/LBMMedium.txt};
        \addplot[only marks,mark=square ,mark size=2pt,color=blue,thick]
        table[x expr={\thisrow{Re}))},y expr=\thisrow{Drag},col sep=comma]{Data/DragCrisis/LBMFine.txt};
        \addplot[only marks,mark=o,mark size=2pt,color=red,ultra thick] 
        table[x expr={\thisrow{Re}))},y expr=\thisrow{Drag},col sep=comma]{Data/Computation.txt};
          \legend{\tiny {Almedij~\cite{almedeij2008drag}},\tiny {Morris ~\cite{morrison2013introduction}},\tiny{Bakic ~\cite{bakic2003experimental}},\tiny{Achenbach~\cite{achenbach1972experiments}},\tiny{LBM Medium~\cite{geier2017parametrization}},
      \tiny{LBM Fine~\cite{geier2017parametrization}},\tiny{Present work}}
        \end{semilogxaxis}
        \end{tikzpicture}
    \end{subfigure}
    \caption{
     Drag crisis: Variation in $C_d$ close to the region of drag crisis. We see a good agreement with the experimental data and past numerical results. \vspace{1 mm}
    } 
    \label{fig: DragSemiLinear}
\end{figure}{}
\begin{figure*}[t!]
\centering
\begin{subfigure}{.32\textwidth}
    \includegraphics[width=1.0\linewidth]{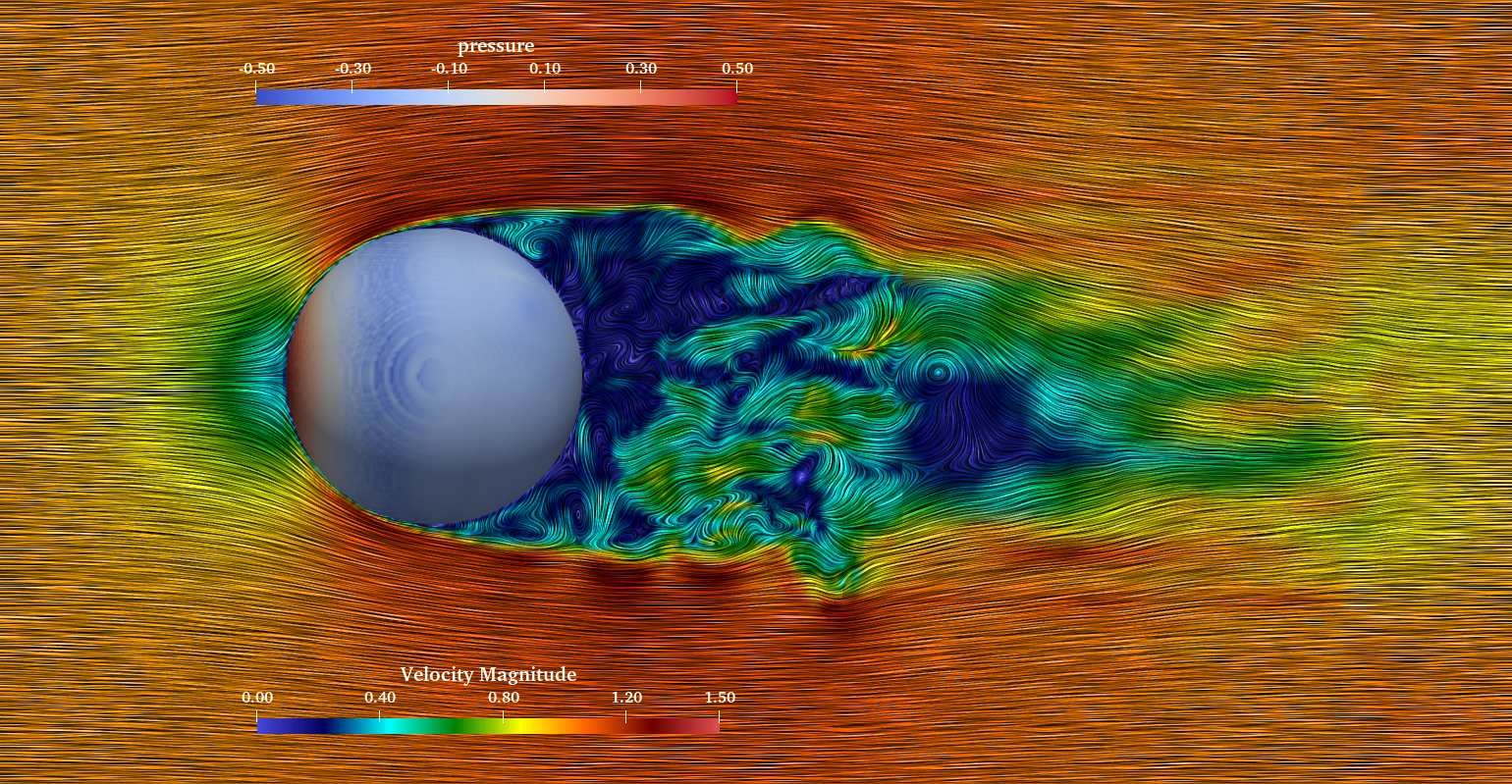}
    \caption{$Re$ = 10,000}
\label{fig: Re10K}
\end{subfigure}
\begin{subfigure}{.32\textwidth}
    \includegraphics[width=1.0\linewidth]{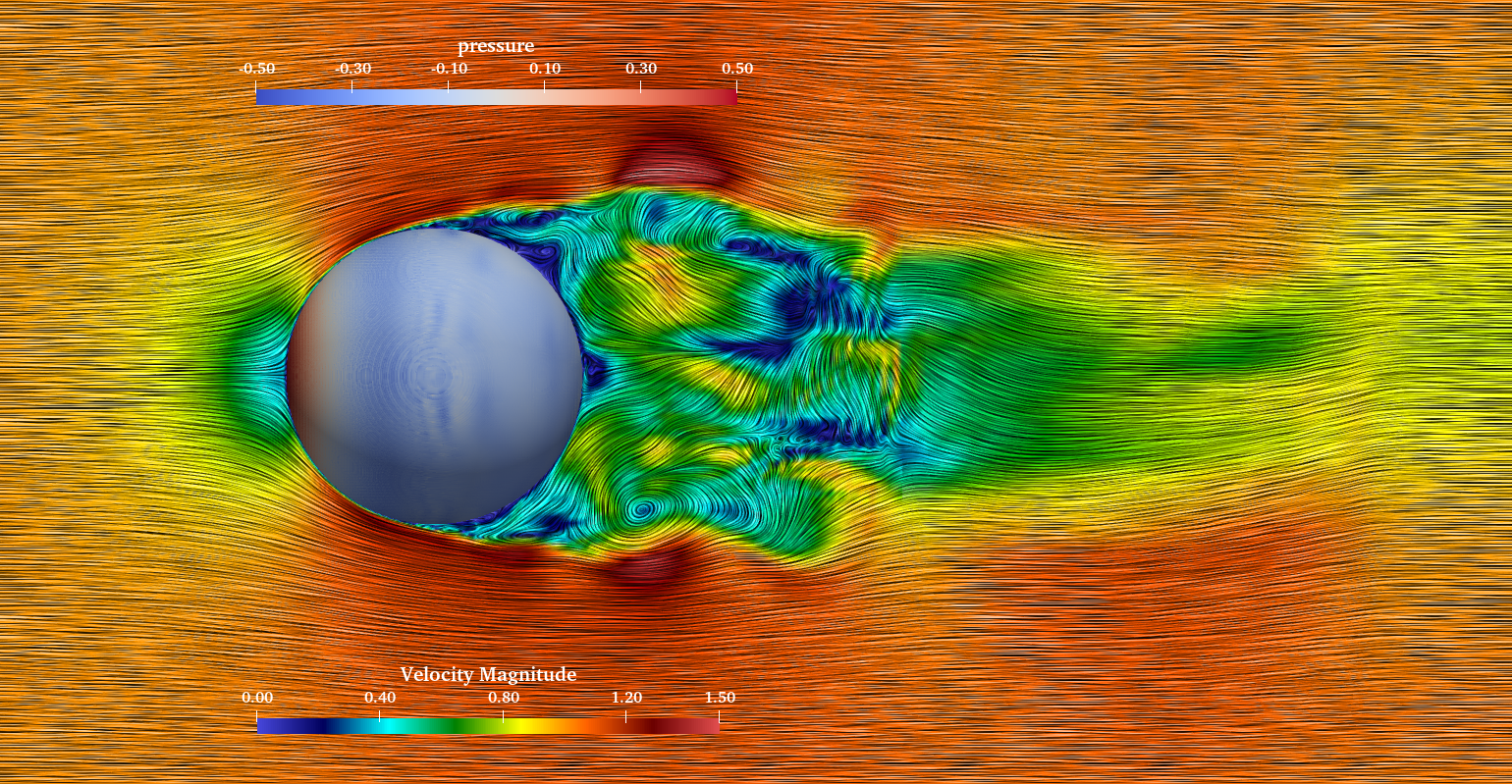}
    \caption{$Re$ = 160,000}
    \label{fig: Re160K}
\end{subfigure}
\begin{subfigure}{.32\textwidth}
    \includegraphics[width=1.0\linewidth]{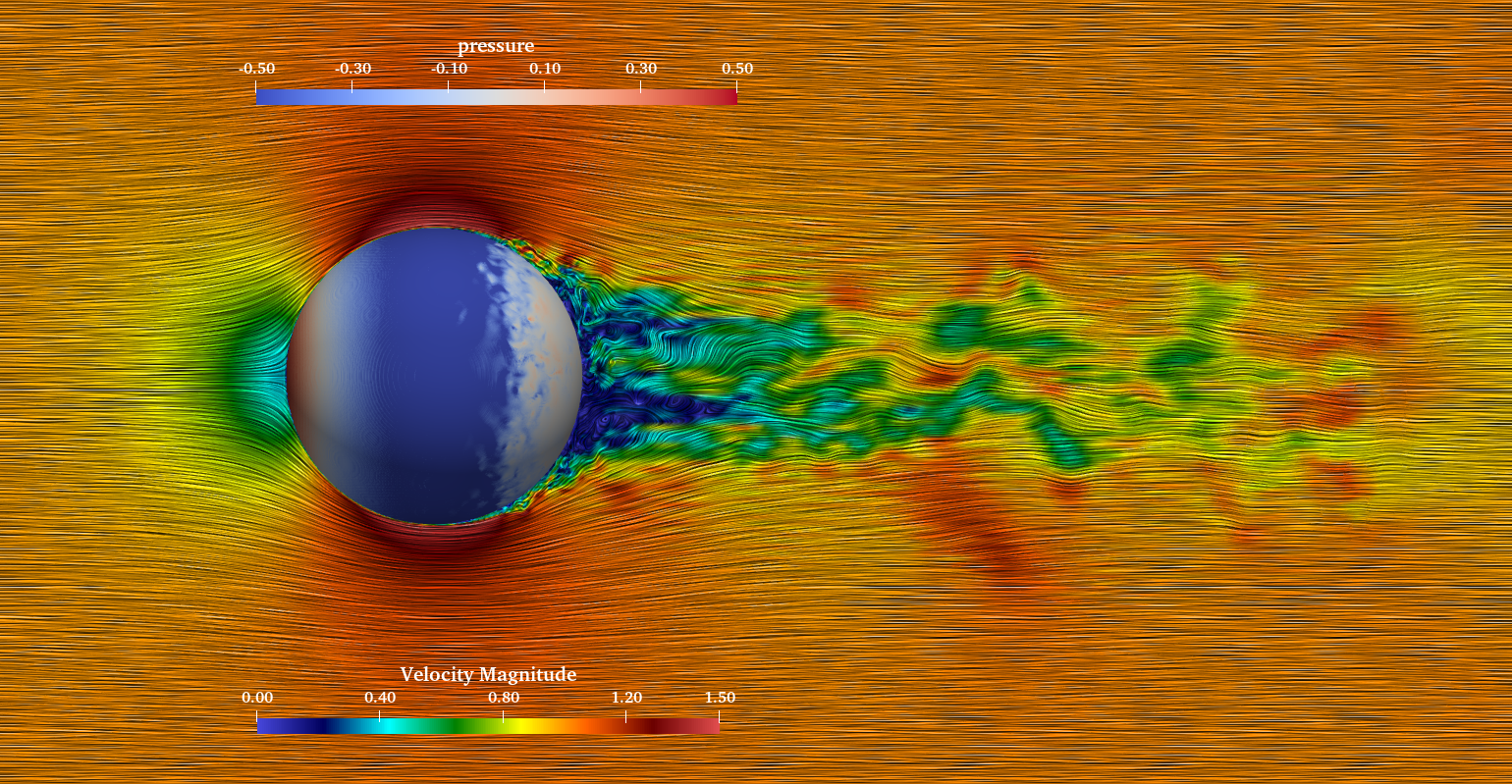}
    \caption{{$Re$ = 2,000,000}}
    \label{fig: Re2M}
\end{subfigure}
\vspace{3mm}
\caption{The wake structures and pressure distribution on sphere at different Reynolds number. The drag crisis is evident by noticing the wake structure as it changes from being divergent at $Re=160,000$ (high drag state) to being convergent at $Re = 2\times10^6$ (low drag state). At the same time, we observe a high pressure region being developed behind the sphere. The development of this high pressure zone is attributed to the low drag state.}
\label{fig: wake}
\end{figure*}
Here, we compare the performance of the proposed algorithm with the existing octree based framework, specifically~\Dendro{}~ \cite{sampath2008dendro,sundar2008bottom,milinda_shayamal_fernando_2020_3876881}. \Dendro~ is a well-validated software and has been widely used in various large scale scientific simulations \cite{mukherjee2010performance,khanwale2020fully,saurabh2020industrial,fernando2018massively,xu2021octree,neilsen2018massively,neilsen2019dendro} and has over 200 citations. \Dendro~ has additional support for carrying out carving operations \cite{xu2021octree}. We choose \Dendro~ as our benchmark for comparison.  For comparison, we choose an elongated channel of dimensions $128 \times 4 \times 1$. The overall mesh is determined by two levels of refinement: base refinement and boundary refinement. Such a channel is commonly found in microfluidic devices, and simulating such devices is an active area of research~\cite{chai2011numerical,shen2004examination}. 

{\renewcommand{\arraystretch}{1.1} 
\begin{table}[t!]
\resizebox{\linewidth}{!}{%
{
\Large
\begin{tabular}{|c|c|c|r|r|r|r|r|}
\hline
\multirow{2}{*}{\begin{tabular}[c]{@{}c@{}}\textbf{Base} \\ \textbf{Refinement}\end{tabular}} & \multirow{2}{*}{\begin{tabular}[c]{@{}c@{}}\textbf{Boundary} \\ \textbf{Refinement}\end{tabular}} & \multirow{2}{*}{\begin{tabular}[c]{@{}c@{}}\textbf{Num} \\ \textbf{Elements}\end{tabular}} & \multirow{2}{*}{\begin{tabular}[c]{@{}c@{}}\textbf{Num}\\ \textbf{Processors}\end{tabular}} & \multicolumn{2}{c|}{\textbf{\Dendro{} (s)}} & \multicolumn{2}{c|}{\textbf{Current Approach (s)}} \\ \cline{5-8} 
 &  &  &  & \begin{tabular}[c]{@{}c@{}}Mesh\\Creation\end{tabular} & \mvec{} & \begin{tabular}[c]{@{}c@{}}Mesh\\Creation\end{tabular} & \mvec{} \\ \hline
\multirow{3}{*}{10} & \multirow{3}{*}{12} & \multirow{3}{*}{3,138,525} & 448 & 107.87 & 59.14 & 1.69 & 9.27 \\ \cline{4-8} 
 &  &  & 896 & 55.41 & 38.03 & 1.00 & 4.53 \\ \cline{4-8} 
 &  &  & 1792 & 38.21 & 25.47 & 0.87 & 2.34 \\ \hline
\multirow{3}{*}{10} & \multirow{3}{*}{14} & \multirow{3}{*}{49,096,209} & 448 & 280.88 & 447.05 & 21.51 & 142.89 \\ \cline{4-8} 
 &  &  & 896 & 159.87 & 349.39 & 11.51 & 77.33 \\ \cline{4-8} 
 &  &  & 1792 & 127.21 & 295.02 & 6.34 & 48.51 \\ \hline
\multirow{3}{*}{12} & \multirow{3}{*}{12} & \multirow{3}{*}{17,440,929} & 448 & - & - & 4.57 & 42.58 \\ \cline{4-8} 
 &  &  & 896 & - & - & 2.29 & 19.65 \\ \cline{4-8} 
 &  &  & 1792 & - & - & 1.75 & 10.76 \\ \hline
\multirow{3}{*}{12} & \multirow{3}{*}{14} & \multirow{3}{*}{63,398,613} & 448 & - & - & 33.56 & 182.23 \\ \cline{4-8} 
 &  &  & 896 & - & - & 17.14 & 125.53 \\ \cline{4-8} 
 &  &  & 1792 & - & - & 9.86 & 47.65 \\ \hline
\end{tabular}
}
}
\vspace{5mm}
\caption{\footnotesize{Comparison of the time (in seconds) for mesh generation and Navier--Stokes \mvec{} for the current approach with \Dendro{} based octree framework. With level $\geq$ 12 of base refinement, ~\Dendro~ framework gave memory error, and hence no time is reported.}}
\label{tab:comparDendro4}
\end{table}
}
Table~\ref{tab:comparDendro4} compares the time for mesh generation and total \mvec{} time\footnote{This time include top-down, bottom up, leaf ~\mvec~ and ghost exchange time} for Navier--Stokes equation. Unlike the 3D Poisson operator used for scaling case, the time to compute elemental operator (denoted by leaf \mvec) is substantially more expensive. The extra overhead introduced by performing top-down and bottom-up traversal  will be significantly smaller compared to performing elemental ~\mvec. The overall time to solve is dominated by load balancing of FEM computation. Since ~\Dendro~ looks at the complete octree, a significant portion of the elements lie inside the void regions. The partitioning algorithm distributes the elements of complete octree (almost) equally among processors. This leads to an imbalance in the overall FEM computation. This is clear from the presented \mvec~ results. Additionally, in all our runs, we were not able to go beyond the level of 12 with the ~\Dendro~ framework. This limits our ability to compare for a further elongated channel. In contrast, within the current framework, we achieve a significant improvement in terms of both scalability and time to solve. \added{Overall, we observed a speedup of about 20$\times$ for mesh generation and 5$\times$ for \mvec~ time. The actual speedup that we can achieve is application-specific and needs to be studied individually. The major factors that determine the overall speedup can be mainly categorized as a) the fraction of volume that can be excluded out; b) complexity of \In--\Out test; c) the resultant communication pattern. }



\section{Application: Classroom Airflow Simulation}\label{sec: ResultApplication}\label{sec: classRoom}

\subsubsection*{\textbf{Validation}:\;} We first validate the solver by demonstrating the ability to capture the drag crisis by simulating flow past a sphere. A sphere of diameter $d=1$ is placed at a distance $3d$ from the inlet at $(3d, \,3d,\, 3d)$ in a computation domain of $(10d, \,6d,\, 6d)$.  The walls of the domain, except the outlet have constant non-dimensional freestream velocity of $(1, \, 0, \, 0)$ and zero pressure gradient. At the outlet, the pressure is set to 0 and zero gradient velocity boundary condition is applied at the wall. No-slip boundary condition (zero Dirichlet) for velocity is imposed on the surface of the sphere. We use a well-established Variational Multiscale (VMS) stabilized Finite Element method for solving the Navier--Stokes equation~\citep{bazilevs2007variational}.

~\figref{fig: DragSemiLinear} plots the variation of $C_d$ across a range of Reynolds numbers close to the drag crisis regime. We see that the results are in excellent agreement with experimental and other numerical results.  We are able to accurately capture the drag crisis phenomena, where a sudden drop in drag from 0.5~\citep{achenbach1972experiments} - 0.6~\citep{bakic2003experimental} at Re around $16,000$ to 0.1~\citep{achenbach1972experiments} - 0.2~\citep{geier2017parametrization} at $Re$ of 2 million is observed. The finest resolution simulation consists of $\sim$40M elements which is significantly lower than for LBM simulation by ~\citet{geier2017parametrization}. We visualize the transition across the drag crises regime in ~\figref{fig: wake}. The drop in drag in ~\figref{fig: DragSemiLinear} is due to the change in wake structure and pressure distribution in ~\figref{fig: wake}.
\begin{figure}[t!]
    \centering
    \includegraphics[trim=0 90 0 100, clip,width=0.99\linewidth]{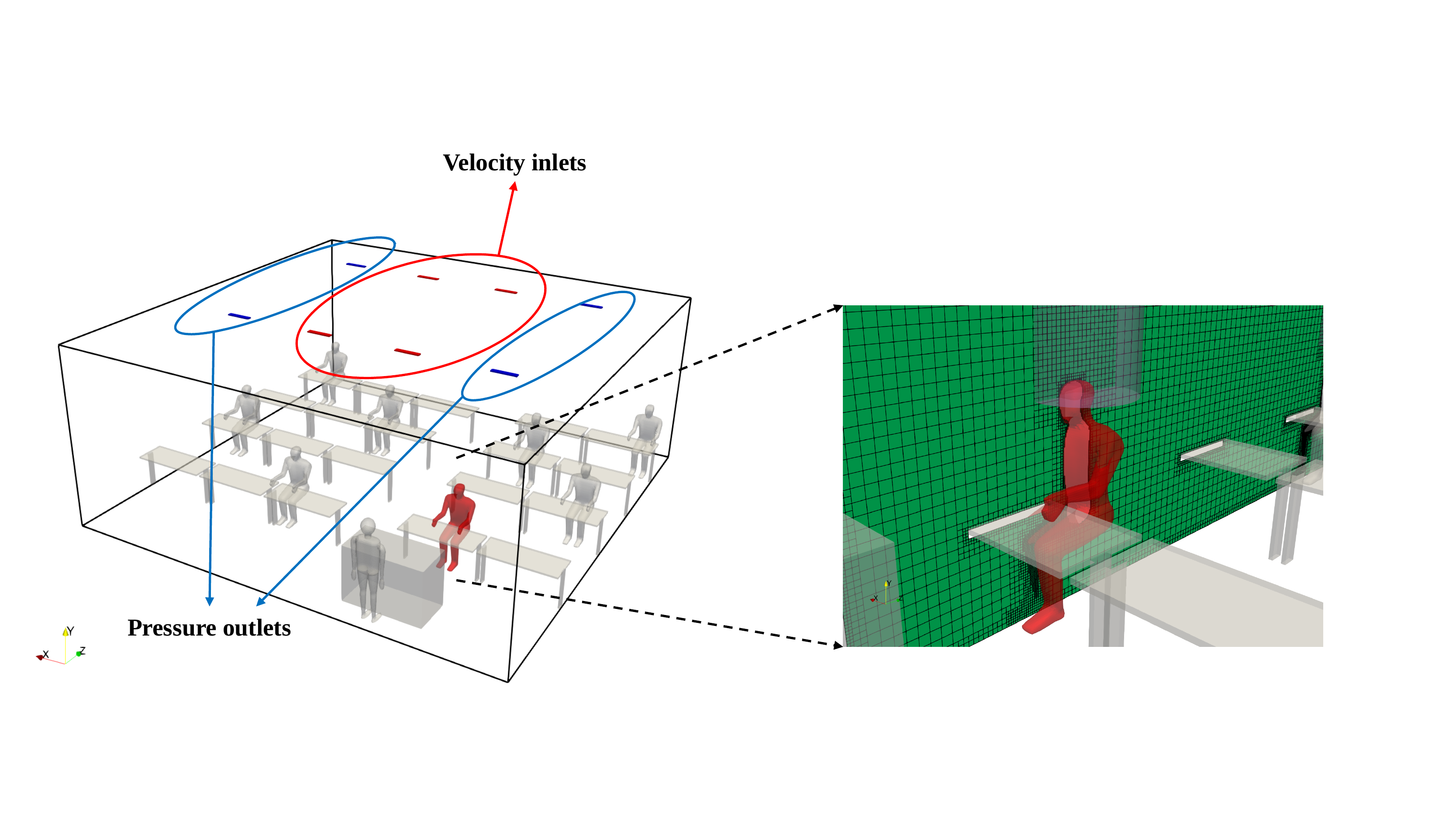}
    \caption{Figure showing the classroom domain and different regions of boundary conditions such as velocity inlet and pressure outlet. The right side shows the zoomed image with mesh refinement near the object. In all our simulations, we consider that the person marked in red is COVID positive.\vspace{2 mm}}
    \label{fig:classroomDomain}
\end{figure}

\begin{figure*}
\centering
\begin{subfigure}{.40\textwidth}
  \centering
  \includegraphics[width=1.0\linewidth,trim=100 180 290 0, clip]{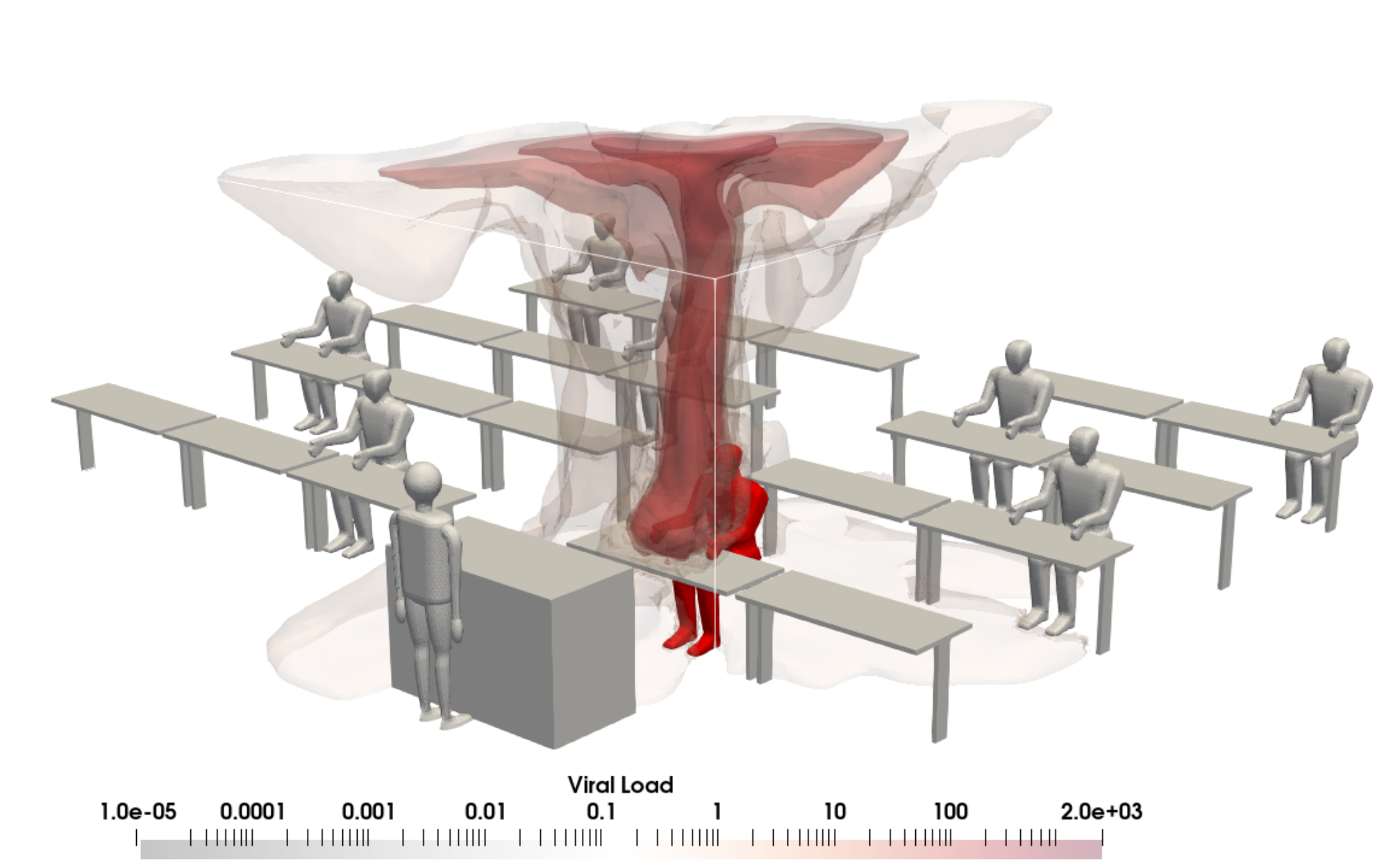}
  \caption{Without monitors}
  \label{fig:classroom:noMonitors}
\end{subfigure}%
\begin{subfigure}{.40\textwidth}
  \centering
  \includegraphics[width=1.0\linewidth,trim=190 180 200 0, clip]{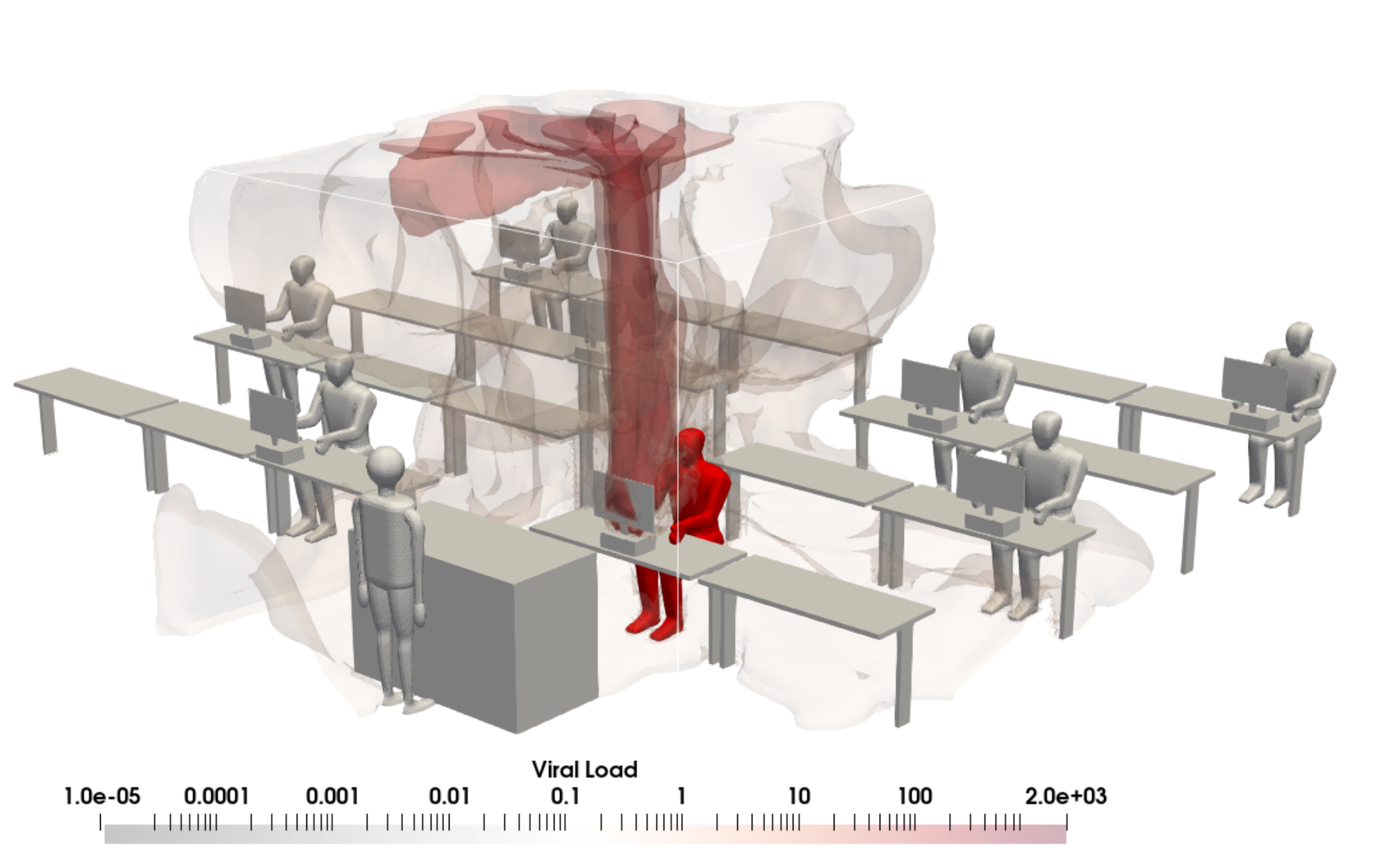}
  \caption{With monitors}
  \label{fig:classroom:Monitors}
\end{subfigure}
\vspace{4 mm}
\caption{\textit{Classroom scenario:} Evaluation of viral load (in quanta / $\mathrm{m}^3$) in two classroom scenarios with and without monitors. The mannequin marked in the red is infected with COVID and transmits the virus. The isocontours represent the regions of different viral load concentrations in space.}
\label{fig:test}
\end{figure*}
\textbf{Application}: We finally demonstrate the ability of our framework to simulate flow past complex geometries. We consider a realistic scenario of modeling airflow in a classroom with complex furniture, seated students with/without computers (and monitors), and a standing instructor. \added{We are particularly interested in accessing if specific locations in the room are at significantly higher risk for transmission -- for example, where there is local recirculation causing limited air exchange with the outside environment. In such cases, it becomes imperative to identify if such locations have a higher risk and rank among alternate seating arrangements to mitigate this risk. The current incomplete octree framework allows us to efficiently and rapidly evaluate various seating arrangements and scenarios. 
} In order to \textit{carve} out the geometry, we perform a series of \In - \Out tests. This gives an automated way to carve out complex geometries from the domain. \figref{fig:classroomDomain} shows the computational domain of size  $4.83 \times 3.34 \times 1$ that includes complex features such as tables, chairs and mannequins representing students and instructor. The velocity inlets and pressure outlets are located at the top of the domain.  The non-dimensional inlet velocity of $(0, 0, -1)$ is imposed at velocity inlets and zero pressure at pressure outlets. $Re = 10^5$ was considered based on the inlet velocity and classroom height. These flow rates, air exchange rates, and room geometry represent typical values in classrooms seen in US schools.

\added{Fig.~\ref{fig:test} illustrates preliminary results enabled by our framework to evaluate the transmission of COVID viral load in the classroom. We evaluate the impact of one infected individual (colored red) who periodically coughs, releasing an aerosolized load of viral particles. We model the time-dependent transmission of the viral load as a scalar transport equation that is advected by a statistically steady--state flow field obtained from the solution of Navier--Stokes solver. We considered a classroom under two different scenarios: with (\figref{fig:classroom:Monitors}) and without (\figref{fig:classroom:noMonitors}) the presence of computer monitors. We observe a significant reduction in transmission risk in the case with monitors due to the monitors redirecting the flow field upwards away from the occupied zone.}
\added{
{\renewcommand{\arraystretch}{1.4} 
\begin{table}[t!]
\resizebox{0.99\linewidth}{!}{%
{
\Large
\begin{tabular}{|c|c|c|c|c|c|c|c|c|c|}
\hline
\multirow{2}{*}{\begin{tabular}[c]{@{}c@{}}\textbf{Base}\\ \textbf{level}\end{tabular}} & \multirow{2}{*}{\begin{tabular}[c]{@{}c@{}}\textbf{Exit} \\ \textbf{refine}\\ \textbf{level}\end{tabular}} & \multirow{2}{*}{\begin{tabular}[c]{@{}c@{}}\textbf{Body}\\ \textbf{refine}\\ \textbf{level}\end{tabular}} & \multicolumn{2}{c|}{\textbf{Elements}} & \multirow{2}{*}{\begin{tabular}[c]{@{}c@{}}\textbf{Num} \\ \textbf{procs}\end{tabular}} & \multicolumn{2}{c|}{\textbf{Immersed (s)}} & \multicolumn{2}{c|}{\textbf{Carved out (s)}} \\ \cline{4-5} \cline{7-10} 
 &  &  & \begin{tabular}[c]{@{}c@{}}\textbf{Active} \\ \textbf{Elements}\end{tabular} & \begin{tabular}[c]{@{}c@{}}$f_{\mathrm{excess}}$ \end{tabular} &  & \begin{tabular}[c]{@{}c@{}}\textbf{Mesh} \\ \textbf{construction}\end{tabular} & \begin{tabular}[c]{@{}c@{}}\textbf{Solve}\\ \textbf{time}\end{tabular} & \begin{tabular}[c]{@{}c@{}}\textbf{Mesh} \\ \textbf{construction}\end{tabular} & \begin{tabular}[c]{@{}c@{}}\textbf{Solve}\\ \textbf{time}\end{tabular} \\ \hline
\multirow{2}{*}{6} & \multirow{2}{*}{8} & \multirow{2}{*}{10} & \multirow{2}{*}{924,549} & \multirow{2}{*}{1.53} & 224 & 92.36 & 178.74 & 38.57 & 61.56 \\ \cline{6-10} 
 &  &  &  &  & 448 & 50.48 & 94.95 & 22.47 & 32.05 \\ \hline
\multirow{2}{*}{6} & \multirow{2}{*}{9} & \multirow{2}{*}{10} & \multirow{2}{*}{1,259,670} & \multirow{2}{*}{1.43} & 224 & 136.13 & 220.6 & 48 & 83.3 \\ \cline{6-10} 
 &  &  &  &  & 448 & 73.24 & 95.06 & 28.88 & 45.10 \\ \hline
\multirow{2}{*}{7} & \multirow{2}{*}{9} & \multirow{2}{*}{11} & \multirow{2}{*}{3,461,548} & \multirow{2}{*}{1.64} & 448 & 210.8 & 309.60 & 89.04 & 131.52 \\ \cline{6-10} 
 &  &  &  &  & 896 & 107.88 & 161.70 & 53.24 & 71.3 \\ \hline
\end{tabular}
}
}
\vspace{0.2in}
\caption{\added{\footnotesize{Comparison of mesh generation and solve time for IBM  with the current approach for the classroom case. The carved out approach leads to a significant reduction in number of elements. $f_{\mathrm{excess}}$ represents the excess fraction of elements obtained as a result of generating complete octree. \vspace{-5 mm}}}}
\label{tab:comparD5}
\end{table}
}
}
\added{
{\renewcommand{\arraystretch}{1.3} 
\begin{table}[t!]
\resizebox{\linewidth}{!}{%
{
\Large
\begin{tabular}{|c|c|c|c|c|c|c|c|c|l|}
\hline
\multirow{2}{*}{\textbf{\begin{tabular}[c]{@{}c@{}}Base \\ level\end{tabular}}} & \multirow{2}{*}{\textbf{\begin{tabular}[c]{@{}c@{}}Exit\\ level\end{tabular}}} & \multirow{2}{*}{\textbf{\begin{tabular}[c]{@{}c@{}}Body \\ Level\end{tabular}}} & \multirow{2}{*}{\textbf{\begin{tabular}[c]{@{}c@{}}Num \\ Elements\end{tabular}}} & \multirow{2}{*}{\textbf{}} & \multicolumn{5}{c|}{\textbf{Number of Processors}} \\ \cline{6-10} 
 &  &  &  &  & 224 & 448 & 896 & 1792 & 3584 \\ \hline
\multirow{2}{*}{7} & \multirow{2}{*}{8} & \multirow{2}{*}{11} & \multirow{2}{*}{5,555,871} & {Time(s)} & 344.87 & 176.15 & 92.16 & 46.44 & 23.95 \\ \cline{5-10} 
 &  &  &  & {Efficiency} & 1.0 & 0.98 & 0.94 & 0.93 & 0.90 \\ \hline
\multirow{2}{*}{9} & \multirow{2}{*}{9} & \multirow{2}{*}{11} & \multirow{2}{*}{23,054,077} & {Time(s)} & - & 539.15 & 272.28 & 142.73 & 75.7 \\ \cline{5-10} 
 &  &  &  & {Efficiency} & - & 1.0 & 0.99 & 0.94 & 0.89 \\ \hline
\end{tabular}
}
}
\vspace{0.2in}
\caption{\added{\footnotesize{\textit{Scaling result for classroom simulation:} Comparison of total solve time and strong scaling efficiency with increase in processor count for two different meshes with varying levels of refinement.}}\vspace{-2 mm}}
\label{tab:scalClassroom}
\end{table}
}
}

\added{\tabref{tab:comparD5} compares the IMGA based immersed implementation (based on the open-source code~\cite{Dendrite})  with the carved-out approach. The overall mesh consists of three refinement levels: base refinement, exit refinement (near the velocity inlet and pressure outlet), and object refinement (near the monitors, tables, and mannequin). To get the refined final mesh, we start the mesh at the base level and successively refine the mesh until the required refinement level is reached. At each iteration, we perform a series of ray-tracing to determine~\In or \Out relative to the object. Overall, we see an approximately 50\% increase in element size for the immersed case.  We achieve a speedup of approximately 2.2$\times$ during the mesh creation stage and 2.8$\times$ for the complete solve time. The speedup obtained in this case is significantly smaller than the channel case described in \secref{sec: comparD4} mainly because of the nature of objects being carved out. The mannequin or the table considered here has a large surface area to volume ratio. The small volume resulted in most of the octants percolating close to the finest level before they can be discarded. Additionally, ray-tracing based \In-\Out test needs to be performed at each iteration of refinement, which is quite expensive \cite{saurabh2020industrial}. Once the mesh is generated, we see a significant speedup in the overall solve time. \tabref{tab:scalClassroom} compares the scaling efficiency for two different meshes. Overall we achieve a good scaling efficiency of about 90\% over 16 fold increase in the number of processors.}
\section{Conclusion} \label{sec: Conclusion}

We present a fast and scalable tree-based mesh generation that is not limited to isotropic domains, which serves as an alternative to using two-tier meshes that are not dependent on having top-level hexahedral meshes. The algorithms presented in the paper are dimension agnostic and provides a generic way to handle any arbitrary geometries.
Our approach allows all elements to remain isotropic, which speeds up assembly and does not affect conditioning due to element stretching. The scaling behavior of the ~\mvec, which is the most dominant part of any FEM solver, has been verified up to $\mathcal{O}(16 K)$ cores. We further showcase the applicability of these algorithms by solving Navier--Stokes for a large-scale 3D problem in the presence of complex geometries. These algorithmic features allows fast, well-balanced creation of complex meshes and efficient solvers that open the way for parametric exploration of very large-scale simulations (as our example simulation suggests). \added{In future, we plan to extend the algorithms to incorporate DG based FEM along with Finite Difference and Finite Volume Methods. }

\newpage
\bibliographystyle{ACM-Reference-Format}
\bibliography{./ms}

\clearpage 

\newpage
\appendix
\section{Artifact Description}
\subsection{Libraries dependencies}
The following dependencies are required to compile the code:
\begin{itemize}
	\item C/C++ compilers with C++11 standards and OpenMP support
	\item MPI implementation (e.g. openmpi, mvapich2 )
	\item \petsc~ 3.8 or higher
	\item ZLib compression library (used to write \texttt{.vtu} files in binary format with compression enabled)
	\item  MKL / LAPACK library
	\item CMake 2.8 or higher version
	\item \href{http://hyperrealm.github.io/libconfig}{libconfig} for parameter reading from file.
\end{itemize}

\subsection{\Frontera~ environment}
Experiments performed in \Frontera~ are executed in the following module environment. 
\begin{verbatim}
Currently Loaded Modulefiles:
  1) intel/19.0.5   4) python3/3.7.0
  2) impi/19.0.5    5) autotools/1.2
  3) petsc/3.12     6) cmake/3.16.1
\end{verbatim}

\subsection{\Frontera~ compute node configuration} \label{sec:Frontera}
\begin{verbatim}
    Architecture:          x86_64
    CPU op-mode(s):        32-bit, 64-bit
    Byte Order:            Little Endian
    CPU(s):                56
    On-line CPU(s) list:   0-55
    Thread(s) per core:    1
    Core(s) per socket:    28
    Socket(s):             2
    NUMA node(s):          2
    Vendor ID:             GenuineIntel
    CPU family:            6
    Model:                 85
    Model name:            Intel(R) Xeon
                           Platinum 8280
                           CPU @ 2.70GHz
    Stepping:              7
    CPU MHz:               2700.000
    BogoMIPS:              5400.00
    Virtualization:        VT-x
    L1d cache:             32K
    L1i cache:             32K
    L2 cache:              1024K
    L3 cache:              39424K
    
    MemTotal:       195920208 kB
    MemFree:        168962328 kB
    MemAvailable:   168337408 kB
\end{verbatim}

\pagebreak
\section{Artifact Evaluation}
\subsection{Signed distance computation}
Let $M$ denotes the closed orientable 2-manifold triangular mesh. The signed distance from a point $\Vector{p}$ rto $M$ is given by:
\begin{equation}
    d(\Vector{p},M) = \inf_{\Vector{x} \in M} \norm{\Vector{p}- \Vector{x}}\mathrm{sign}(\Vector{n}\cdot(\Vector{p} - \Vector{c}))
\end{equation}
where: $\Vector{c}$ denotes the closest point to $\Vector{p}$ and $\Vector{n}$ denotes the outward normal. A positive value of $d$, means the point is inside the surface and vice-versa.
\subsection{Details of solver selection}
\petsc~ was used to solve all the linear algebra problems. In particular, bi-conjugate gradient descent (\texttt{-ksp\_type bcgs}) solver was used in conjunction with Additive - Schwartz (\texttt{-pc\_type asm}) preconditioner to solve the linear system of equations. The \textsc{NEWTONLS} class by \petsc, that implements a Newton Line Search method, was used for the nonlinear problems. Both the relative residual tolerance and the absolute residual tolerance for linear and non - linear solve are set to $10^{-6}$ in all numerical results.







\subsection{Downloading and installing the code}
This section presents how to run and reproduce the results presented in the paper. You can clone the repository using, \texttt{git clone bitbucket.org/baskargroup/sc21-kt.git}. We use \texttt{CMake} to configure and build. In \Frontera~node, 
\begin{itemize}
    \item \texttt{git clone bitbucket.org/baskargroup/sc21-kt.git}
    \item \texttt{mkdir build \&\& cd build}
    \item Load the module environment
    \item \texttt{mkdir build \&\& cd build}
    \item \texttt{cmake ../.}
    \item \texttt{make MVCChannel MVCSphere signedDistance}
\end{itemize}

\subsection{Running experiments}
\subsubsection{\texttt{MVCChannel}: \;} Scaling  run for channel case. In order to run the scaling case, it requires the 3 parameters: a) baseLevel b) boundaryLevel and c) element order (1/2). For example, on ~\Frontera~ it can be ran as:

\texttt{ibrun MVCChannel 10 12 1 log10\_12.out} 

\noindent where, 10 is the base refinement level, 12 is the boundary refinement level, 1 is the element order and \texttt{log10\_12.out} is the output file containing the relevant timing information.

\subsubsection{\texttt{MVCSphere}: \;}
Scaling  run for sphere case. In order to run the scaling case, it requires the 3 parameters: a) baseLevel b) boundaryLevel and c) element order (1/2). For example, on ~\Frontera~ it can be ran as:

\texttt{ibrun MVCSphere 7 12 1 log7\_12.out} 

\noindent where, 7 is the base refinement level, 12 is the boundary refinement level, 1 is the element order and \texttt{log7\_12.out} is the output file containing the relevant timing information.

\subsubsection{\texttt{signedDistance}: \;} The computation of signedDistance. In order to find the signed distance, we successively refined near the boundaries and computed the signed distance. In order to run the code:

\texttt{ibrun signedDistance stlFileName 4 14}

\noindent where: \texttt{stlFileName} is the name of stl file, 4 is the minimum refinement level and 14 is the maximum level of refinement at the stl boundary. After each successive iteration, the code output the information of boundary nodes. Then we compute the signed distance, using the python script provided under the \texttt{scripts} folder. In order the run the python scripts:

\texttt{python3 signedDistance stlFileName}.

Note that you might need to change the number of processor on your machine as the python script is parallel and make use of \texttt{multiprocessing} library.

\subsubsection{Roofline plot: \:}

We computed the roofline using Intel Advisor. In order to run the roofline plot, first run the survey using:

\noindent\texttt{ibrun -np 1 advixe-cl -collect survey -project-dir outputDir
           -search-dir \\
           src:=examples/BenchMark\_channel/src 
           -- MVCChannel \\
           baseLevel boundaryLevel eleOrder outputFile}
           
\noindent where: 
\text{outputDir} is the directory to store output and \texttt{MVCChannel   baseLevel boundaryLevel eleOrder outputFile} is the same as in previous description of Channel scaling.

Finally in order to collect FLOPS count:

\noindent\texttt{ibrun -np 1 advixe-cl -collect=tripcounts
           --flop  \\
            --mark-up-list= src/benchmark.cpp 
            -project-dir= outputDir 
            -- MVCChannel baseLevel \\
            boundaryLevel eleOrder outputFile}
            
\noindent where:
\text{outputDir} must be same directory as the survey directory and \texttt{MVCChannel} must be called with the same arguments.

\end{document}